\documentclass[10pt]{amsart}
\usepackage[utf8]{inputenc}
\usepackage{geometry}
\usepackage{amsmath,amssymb,amsthm,mathtools}
\usepackage{natbib}
\usepackage{xcolor}
\usepackage{booktabs}
\usepackage{float}
\usepackage{colortbl}
\usepackage{hyperref}

\theoremstyle{definition}
\newtheorem{theorem}{Theorem}[section]
\newtheorem{proposition}[theorem]{Proposition}
\newtheorem{lemma}[theorem]{Lemma}
\newtheorem{corollary}[theorem]{Corollary}
\newtheorem{remark}[theorem]{Remark}
\numberwithin{equation}{section}

\newcommand{\kt}[1]{\textup{(#1)}}                 
\newcommand{\QQ}{\mathbb{Q}}
\newcommand{\ZZ}{\mathbb{Z}}
\newcommand{\CC}{\mathbb{C}}
\newcommand{\RR}{\mathbb{R}}
\newcommand{\PP}{\mathbb{P}}
\newcommand{\NN}{\mathbb{N}}
\newcommand{\cO}{\mathcal{O}}

\newcommand{\cI}{\mathcal{I}}
\newcommand{\cM}{\mathcal{M}}
\newcommand{\Mbar}{\overline{\mathcal{M}}}

\newcommand{\bone}{\mathbf{1}}
\newcommand{\ceq}{\coloneqq}
\DeclareMathOperator{\Gr}{Gr}

\DeclareMathOperator{\rk}{rk}

\DeclareMathOperator{\id}{id}
\DeclareMathOperator{\End}{End}
\DeclareMathOperator{\NS}{NS}

\DeclareMathOperator{\NE}{NE}

\DeclareMathOperator{\ev}{ev}
\newcommand{\tors}{/\mathrm{tors}}
\newcommand{\Hdg}{\mathrm{Hdg}}

\newcommand{\pt}{\mathrm{pt}}

\newcommand{\odd}{\mathrm{odd}}
\newcommand{\even}{\mathrm{even}}

\begin{document}
\address{Department of Mathematics, Imperial College London,
  London SW7 2AZ, United Kingdom}
\email{aideen.fay23@imperial.ac.uk}
\title{Quantum Cohomology and irrationality of K\"uchle fourfolds}
\author{Aideen Fay}
\date{\today}

\begin{abstract} 

    We prove that very general K\"uchle fourfolds of types \kt{c1}, \kt{c6}, \kt{b8}, and \kt{b11} are irrational. The key ingredients are a new rationality criterion, based on the theory of \emph{Hodge atoms} due to Katzarkov--Kontsevich--Pantev--Yu and its integral extension by Gu\'er\'e, and a computation of the action of quantum multiplication on the vanishing cohomology of these fourfolds, which uses a version of the Quantum Lefschetz Hyperplane Principle. Quantum multiplication acts as a scalar on the vanishing cohomology, and we give an enumerative interpretation of this scalar.
\end{abstract}

\maketitle

\section{Introduction}

In~\cite{kuchle}, K\"uchle classified all smooth index-one Fano fourfolds that arise as the zero locus of a regular section of a globally generated, completely reducible homogeneous vector bundle on a Grassmannian. 
Our main result concerns their rationality. 
\begin{theorem}\label{cor:kuchleirr}
  Very general K\"uchle fourfolds of types \kt{b8}, \kt{b11}, \kt{c1}, and \kt{c6} are irrational. 
\end{theorem}
What makes the birational geometry of these fourfolds particularly compelling is that their middle cohomology carries a large surface-type Hodge structure that is not of K3 type. Whereas the obstructions to rationality for cubic, Gushel--Mukai, and Verra fourfolds~\cite{Gueirrational,
benedetti2026quantumcohomologyirrationalitygushelmukai,Gueintegral} involve K3 surfaces, rationality of a very general fourfold of type~\kt{c1} or ~\kt{c6} would imply the existence of a projective surface $S$ of geometric genus $5$ or $2$, respectively.

The proof of Theorem~\ref{cor:kuchleirr} rests on a rationality criterion for Fano fourfolds, which we now state. 
Let $X$ be a smooth projective complex fourfold and define
$$
\Lambda_X\ceq H^4(X,\ZZ)\tors,
\qquad
A_X\ceq\Lambda_X\cap H^{2,2}(X),
\qquad
T_X\ceq A_X^\perp\subseteq\Lambda_X,
$$
where orthogonality is taken with respect to the intersection pairing $(\alpha,\beta)_X=\int_X\alpha\smile\beta$; we write $T_{X,\QQ}\ceq T_X\otimes_\ZZ \QQ$. For a smooth projective surface $S$, write
$$
\Lambda_S\ceq H^2(S,\ZZ)\tors,
\qquad
T_S\ceq(\NS(S)\tors)^\perp\subset\Lambda_S,
\qquad
\rho(S)\ceq\rk(\NS(S)\tors)
$$
We use the following notation for a nondegenerate lattice $L$, as recalled in Section~\ref{sec:lattices}:
\begin{align*}
    & D(L) && \text{the discriminant group $L^\vee/L$}\\
    & \ell(D(L)) && \text{the minimal number of generators of $D(L)$}
\end{align*}
As we will describe in Section~\ref{subsec:ql} below, the cohomology group $H^\bullet(X;K)$, where $K$ is a suitable coefficient field, carries an endomorphism $U_X(\tau)$ given by quantum multiplication by the Euler vector field, depending on parameters $\tau$ in the direction of rational Hodge classes. If $X$ satisfies \eqref{eq:conditions} below, then $U_X(\tau)$ acts on the transcendental cohomology $T_X \otimes_{\ZZ} K$ by a scalar, as we show in the proof of Proposition~\ref{pro:quantum-integral}: let us call this scalar $\lambda_X(\tau)$.  Put $\lambda_X^0=\lambda_X(0)$, and define, for general $\tau$,
$$
m_X\ceq\operatorname{mult}_{\lambda_X(\tau)}
       \det\bigl(t\,\id-U_X(\tau)\bigr),
\qquad
M_X\ceq\operatorname{mult}_{\lambda_X^0}
       \det\bigl(t\,\id-U_X(0)\bigr).
$$
Here $m_X$ is the dimension of the generalized eigenspace of $U_X(\tau)$ containing $T_X$, while $M_X$ is the specialization of this quantity to the small quantum locus $\tau=0$. Since multiplicities can only jump up at special points, $m_X\le M_X$.


\begin{proposition}[Quantum Noether integral obstruction]
\label{pro:quantum-integral}
Let $X$ be a smooth projective complex fourfold satisfying
\begin{equation}
\label{eq:conditions}
h^{2,0}(X)=h^{4,0}(X)=0,\qquad h^{3,1}(X)>0,\qquad
\End_{\mathrm{Hdg}}(T_{X,\QQ})
=\QQ\,\mathrm{id}.
\end{equation}
Let $r=\rk T_X$. If $X$ is rational, there is a smooth minimal projective surface $S$
and an integral Hodge isometry $T_X(1)\simeq T_S$ such that:
\begin{enumerate}
\item We have
$$\dim_{\QQ}H^\bullet(S,\QQ)\leq m_X\leq M_X$$
and in particular,
$
r+\max\{1,\ell(D(T_X))\}+2\leq m_X\leq M_X.
$
\item The Picard-rank $\rho(S)$ of $S$ satisfies
$$
\ell(D(T_X))\leq\rho(S)\leq 10+12h^{3,1}(X)-r.
$$
\item If $M_X\leq r+3$, then
$$
\rho(S)=1,\qquad h^{1,0}(S)=0,\qquad p_g(S)=h^{3,1}(X),
$$
the group $D(T_X)$ is cyclic, and
\begin{equation}
\label{eq:canonical-square}
|D(T_X)|k^2=9+12h^{3,1}(X)-r
\qquad\text{for some }k\in\ZZ.
\end{equation}
\end{enumerate}
\end{proposition}

The criterion also recovers the recent result of Gu\'er\'e: 
\begin{corollary}[Gu\'er\'e, {\cite[Corollary~28]{Gueintegral}}]\label{cor:verra}
A very general Verra fourfold is irrational.
\end{corollary}

\begin{table}[htbp]
    \centering
    \renewcommand{\arraystretch}{1.15}
    \setlength{\tabcolsep}{6pt}
    \begin{tabular}{@{}cccccc@{}}
        \toprule
        K\"uchle Type
        & $h^{3,1}(X)$
        & $\rk T_X$
        & $M_X$
        & $D(T_X)$
        & Rationality Criterion \\
        \midrule
        \kt{b8}
        & $3$ & $42$ & $45$
        & $\ZZ/2\oplus\ZZ/6$
        & $2\le \rho(S)\le 1$ \\

        \kt{b11}
        & $2$ & $32$ & $35$
        & $(\ZZ/3)^2$
        & $2\le \rho(S)\le 1$ \\

        \kt{c1}
        & $5$ & $60$ & $63$
        & $\ZZ/21$
        & $21k^2=9$ \\

        \kt{c6}
        & $2$ & $32$ & $35$
        & $\ZZ/9$
        & $9k^2=1$ \\

        \bottomrule\\
    \end{tabular}
\caption{Rationality obstructions for very general K\"uchle fourfolds obtained from Proposition~\ref{pro:quantum-integral}.}
\end{table}

\subsection*{Related work}

Amongst the twenty-one families of K\"uchle varieties~\cite{kuchle}, only those of types \kt{c5}, \kt{c7}, and \kt{d3} are of K3 type~\cite[Section~2.3]{AIF_0__0_0_A50_0}. Kuznetsov studied the geometry of the varieties of types \kt{b4}, \kt{b9}, \kt{c7}, and \kt{d3} in~\cite{Kuznetsov_2015}, where he proved type \kt{b9} to be rational \cite[Theorem~2.5]{Kuznetsov_2015}.  

Kuznetsov also studied the birational geometry of the K\"uchle fivefold of type \kt{c5} in \cite{Kuznetsov_2016}. Taking its smooth hyperplane section gives the K\"uchle fourfold of type \kt{c5}, whose very general member was proven to be irrational by Benedetti--Fay--Gu\'er\'e--Manivel--Perrin~\cite[Proposition~5.1]{benedetti2026atomiccriterionirrationalityquantum}. This result uses Noether's formula to bound the rank of
the transcendental part of a surface atom \cite[Theorem~4.1]{benedetti2026atomiccriterionirrationalityquantum}; here we combine it with a dimension bound from small quantum multiplication and integral lattice data. The type \kt{c5} irrationality result was recovered by Alexeev--Cavenaghi--Galindo--Georgiev--Katzarkov--Martins in~\cite{alexeev2026atomicdecompositioncompleteintersection}, who also prove irrationality of very general K\"uchle fourfolds of types \kt{c2} and \kt{b10} \cite[Theorem~6.13 and 6.16]{alexeev2026atomicdecompositioncompleteintersection}. Both papers have related monodromy arguments concerning the scalar action on vanishing cohomolology; see \cite[proof of~Theorem~4.1]{benedetti2026atomiccriterionirrationalityquantum} and \cite[Theorem~5.35]
{alexeev2026atomicdecompositioncompleteintersection}.

For the fourfold of type $\kt{c7}$, Kuznetsov identified it with the blowup of a cubic fourfold along a Veronese surface in~\cite[Corollary~4.11]{Kuznetsov_2015}. Using this description, the author proved a very general K\"uchle fourfold of type \kt{c7} to be irrational~\cite[Corollary~1.4]{fay2026hasseobstructionsrationalityspecial}. 

The motivating thread throughout this paper's work on birational geometry is the work of Katzarkov--Kontsevich--Pantev--Yu in~\cite{katzarkov2026birationalinvariantshodgestructures}. To apply it to the fourfolds considered here, we use Gu\'er\'e's recent extension to integral Hodge structures~\cite[Proposition~26]{Gueintegral}. 
Our approach for understanding and recovering the relevant quantum cohomology is driven by the Quantum Lefschetz and Riemann--Roch methods of Coates--Givental~\cite{coates2014quantumlefschetzprinciplevector, Coates2001QuantumRL}. Our motivation for applying it within the context of Fano fourfolds finds itself in the large number of quantum periods computed by Coates--Galkin--Kasprzyk--Strangeway in~\cite{Coates01062020}.

\subsection*{Outline}
In Section~\ref{sec:monodromy} we show that quantum multiplication by the hyperplane class $H$ acts by a scalar on vanishing cohomology, and that this scalar is a repeated eigenvalue on the inherited part of cohomology with a nonzero isotropic eigenvector. In Section~\ref{sec:lines}, we express this scalar as a line count and construct embeddings of vanishing and transcendental cohomology using the Fano surface of lines. We recall the Quantum Lefschetz theorem in Section~\ref{sec:comp} and apply it to the K\"uchle fourfolds of type \kt{c1}, \kt{c6}, \kt{b8}, and \kt{b11}. In Section~\ref{sec:lattices} we prove our rationality criterion, which we apply to these K\"uchle fourfolds in Section~\ref{sec:results}. We also apply it to Verra fourfolds in Section~\ref{sec:verra}.
\medskip
All varieties throughout this paper are smooth and projective over $\CC$.

\section{Monodromy and Vanishing Cohomology}\label{sec:monodromy}
Let $X\subset Y\subset\mathbb P^N$ be a smooth hyperplane section of a smooth projective fivefold, and let $i\,\colon X\hookrightarrow Y$ denote the inclusion. Assume that $X$ is an index-one Fano fourfold of Picard-rank one with $H^{\mathrm{odd}}(X,\QQ)=0$. Write $H=c_1(X)$ and set $K=\QQ(q)$.  

Define the inherited and vanishing cohomology by
$$
B_X\ceq\operatorname{im}\!\left(
i^*\colon H^\bullet(Y,\QQ)\longrightarrow
H^\bullet(X,\QQ)\right),
\qquad
V_X\ceq\left(B_X\cap H^4(X,\QQ)\right)^\perp
\subset H^4(X,\QQ),
$$
where orthogonal complements are taken with respect to the Poincar\'e pairing. By weak and hard Lefschetz, the pairing on $B_X$ is nondegenerate and we have the orthogonal decomposition 
$$
H^\bullet(X,\QQ) =H^{\mathrm{even}}(X,\QQ) = B_X \oplus V_X.
$$

Let $U\subset|\mathcal O_Y(1)|$ be the locus of smooth hyperplane sections and $u_0 = [X]\in U$. The local system of vanishing cohomology over $U$ defines the monodromy representation
$$
\rho_{\mathrm{van}}\colon \pi_1(U,u_0) \longrightarrow O(V_X),
\qquad
\Gamma \ceq \operatorname{im}(\rho_{\mathrm{van}}),
$$
where $\Gamma$ is the \emph{monodromy group}; see \cite[\S 5]{deligne}.

\subsection*{Gromov--Witten invariants}
We briefly recall stable maps and Gromov--Witten invariants; see \cite{Kontsevich_1994, li1998virtualmodulicyclesgromovwitten}.  For a smooth projective variety $Z$ and a class $\beta\in H_2(Z,\ZZ)$, let $\overline{\mathcal M}_{0,n}(Z,\beta)$ be the moduli space of \emph{stable maps} 
$$f\colon (C,p_1,\dots,p_n)\to Z\quad \text{with}\qquad f_*[C]=\beta.$$ 
Here $C$ is a connected projective nodal curve of arithmetic genus-zero with distinct smooth marked points $p_1,\dots,p_n\in C$, and every irreducible component of $C$ contracted by $f$ contains at least three special points (marked points or nodes).

For each $j$, there is an \emph{evaluation map}
\[
 \ev_j\colon \overline{\mathcal M}_{0,n}(Z,\beta) \rightarrow Z, \qquad \ev_j(f) = f(p_j).
\]
The space $\overline{\mathcal M}_{0,n}(Z,\beta)$ carries a virtual fundamental class $[\overline{\mathcal M}_{0,n}(Z,\beta)]^{\mathrm{vir}}$ of expected dimension $\dim Z-3+n+c_1(Z)\cdot\beta$. 

The genus-zero \emph{Gromov--Witten invariants} of $Z$ are then the rational virtual intersection numbers
\[
 \bigl\langle\alpha_1\psi_1^{k_1},\dots,\alpha_n\psi_n^{k_n}\bigr\rangle^{Z}_{0,n,\beta}\ceq\int_{[\overline{\mathcal M}_{0,n}(Z,\beta)]^{\mathrm{vir}}}\,\prod_{j=1}^{n}\ev_j^*\alpha_j\smile\psi_j^{k_j}, \qquad\alpha_j\in H^{\bullet}(Z,\QQ),
\]
where $\psi_j$ is the first Chern class of the line bundle whose fibre at $f$ is the cotangent line of $C$ at $p_j$; until Section~\ref{sec:comp} all $k_j$ are zero. Note that these invariants do not change when $Z$ is deformed, and they satisfy the \emph{divisor equation}
\[
\langle D,\alpha_1,\dots,\alpha_n\rangle^{Z}_{0,n+1,\beta}=(D\cdot\beta)\,\langle\alpha_1,\dots,\alpha_n\rangle^{Z}_{0,n,\beta}
\]
for $D\in H^2(Z,\QQ)$ and $\beta\neq0$. When $Z$ carries a hyperplane class $H$, we write an integer $d$ in place of $\beta$ for the disjoint union of moduli spaces, or the sum of the invariants, over all $\beta\in H_2(Z,\ZZ)$ with $H\cdot\beta=d$.

\subsection*{The quantum product}
For the fourfold $X$ with $H=c_1(X)$, the \emph{small quantum product} on $H^\bullet(X,K)$ is the $K$-billinear product defined by 
\[
  (a\star b,c)=\sum_{d\ge0}q^d\,\langle a,b,c\rangle^{X}_{0,3,d}\qquad\text{for all } a,b,c\in H^{\bullet}(X,\QQ).
\]
Here the Novikov variable $q$ records the degree and the term with $d=0$ is $\int_Xa\smile b\smile c$. By the dimension formula, the coefficient of $q^d$ in $a\star b$ has degree $\deg a+\deg b-2d$, so the sum is finite. The product $\star$ is graded commutative, associative, and has unit $1 \in H^0(X,K)$. 

Write
$$
B\ceq B_X\otimes_{\mathbb{Q}}K,\qquad V\ceq V_X\otimes_{\mathbb{Q}}K.
$$
The small quantum product $\star$ makes
$$
R\ceq H^\bullet(X,K)=B\oplus V
$$
a \emph{Frobenius algebra}: that is, a finite-dimensional unital commutative associative $K$-algebra equipped with a nondegenerate symmetric bilinear pairing satisfying $(a \star b, c) = (a, b \star c)$.   

This compatibility, together with the monodromy equivariance of $\star$ proved in Proposition~\ref{prop:orthogonal-monodromy}, leads to Proposition~\ref{prop:character}: an abstract Frobenius algebra version of Hu's eigenvector construction for projective
complete intersections; see \cite[Proposition~5.1]{hu2025bigquantumcohomologyfano}. 
\begin{proposition}\label{prop:orthogonal-monodromy}
If the pairing on $V_X\otimes_\QQ \RR$ is indefinite, then $\Gamma$ is Zariski dense in the orthogonal group $O(V_X)$. Consequently, the small quantum product on $R$ is $O(V)$-equivariant---with $O(V)$ acting trivially on $B$ and by its  standard representation on $V$.
\end{proposition}
\begin{proof}
By standard theory, $\Gamma$ is generated by Picard--Lefschetz reflections in the vanishing cycles of a Lefschetz pencil through $X$. These cycles span $V_X$ and form a single orbit up to sign. Thus  $V_X$ is an absolutely irreducible $\Gamma$-module. 

By Deligne's reflection-group alternative \cite[Lemme~4.4.2s]{deligne1980weil2} (see also \cite[Part~II]{katz2004larsen} for a statement in English), the Zariski closure of $\Gamma$ is either a finite reflection group or $O(V_X)$. It cannot be finite, as assuming so contradicts the stated indefiniteness: a finite group preserves a positive definite form on $V_X\otimes\mathbb{R}$, which by absolute irreducibility and Schur's Lemma is proportional to the pairing, which would then be definite. 

Finally, deformation invariance of Gromov--Witten invariants makes the small quantum product equivariant under $\Gamma$. This is a Zariski-closed condition, and hence extends to $O(V_X)$ and, after extending scalars, to $O(V)$. 
\end{proof}

\begin{proposition}\label{prop:character}
Let $R=B\oplus V$ be an orthogonal decomposition of a Frobenius algebra over a field $k$ of characteristic zero, with $1\in B$ and $\dim_k V\ge2$, whose product $\star$ is $O(V)$-equivariant as in Proposition~\ref{prop:orthogonal-monodromy}. 

Then $B$ is a subalgebra, and there are a unique character $\chi\colon B\to k$ and a unique class $b_\chi\in B$ such that
\begin{equation}\label{eq:chi}
a\star v=\chi(a)\,v,\qquad v\star w=(v,w)\,b_\chi,\qquad (b_\chi,a)=\chi(a)
\end{equation}
for all $a\in B$ and $v,w\in V$; in particular $(b_\chi,1)=1$. Moreover,
\begin{equation}\label{eq:bchi}
a\star b_\chi=\chi(a)\,b_\chi\qquad\text{and}\qquad \chi(b_\chi)=0,
\end{equation}
so that $b_\chi$ is isotropic and $b_\chi\star b_\chi=0$. Finally, for $H\in B$, writing $L_B\ceq(H\star)|_B$ and $\lambda\ceq\chi(H)$,
\begin{equation}\label{eq:double-root}
(t-\lambda)^2\ \text{ divides }\ \det\bigl(t\,\id_B-L_B\bigr).
\end{equation}
\end{proposition}

\begin{proof}
Equivariance of the quantum product under $\mathrm{id}_B\oplus(-\mathrm{id}_V)$ gives
$$
B\star B\subseteq B,\qquad
B\star V\subseteq V,\qquad
V\star V\subseteq B.
$$
The standard representation of  $O(V)$ is absolutely irreducible, so Schur's lemma gives $a\star v=\chi(a)v$. Since $\star$ is associative and unital, $\chi$ is an algebra homomorphism. Nondegeneracy of the pairing on $B$ provides a unique $b_\chi\in B$ with $(b_\chi,a)=\chi(a)$ for all $a\in B$. For $v,w\in V$ and $a\in B$, the Frobenius property then gives
$$
(v\star w,a)=(v,w\star a)=\chi(a)\,(v,w)=\bigl((v,w)\,b_\chi,\,a\bigr),
$$
so $v\star w=(v,w)\,b_\chi$, again by nondegeneracy. This proves \eqref{eq:chi}.

Now fix $v,w\in V$ with $(v,w)=1$, so that $b_\chi=v\star w$. For $a\in B$, associativity gives $a\star b_\chi=(a\star v)\star w=\chi(a)\,b_\chi$. Since $\dim V\ge 2$, there is a nonzero $u\in V$ orthogonal to $w$; then $w\star u=0$, so
$$
\chi(b_\chi)\,u=b_\chi\star u=(v\star w)\star u=v\star(w\star u)=0
$$
and therefore $\chi(b_\chi)=0$. This proves \eqref{eq:bchi}, whence $(b_\chi,b_\chi)=\chi(b_\chi)=0$ and $b_\chi\star b_\chi=\chi(b_\chi)\,b_\chi=0$.

Finally, $b_\chi\ne 0$ because $(b_\chi,1)=\chi(1)=1$, so $k\,b_\chi$ is a line in the hyperplane $\ker\chi$. Both are preserved by $L_B$, which acts by $\lambda$ on $k\,b_\chi$ and on $B/\ker\chi$, since $H\star b_\chi=\lambda\,b_\chi$ and $\chi(H\star a)=\lambda\,\chi(a)$ for $a\in B$. This gives \eqref{eq:double-root}.
\end{proof}

When $\det(t-L_B)$ has more than one repeated root, as for types \kt{b11} and \kt{c1} below, \eqref{eq:chi} and \eqref{eq:bchi} decide between them, as the eigenspace of $\lambda$ contains the isotropic class $b_{\chi}$. For an example of this, see the proof of Proposition~\ref{prop:quantum-b11}.

\section{Lines and the Cylinder Correspondence}\label{sec:lines}
Keeping the notation of the previous section, we identify the scalar $\lambda$ as $-q$ times the number of lines of $Y$ through a general point. Although we do not use this result in the proof of Theorem~\ref{cor:kuchleirr}, where $\lambda$ is found from \eqref{eq:bchi} and \eqref{eq:double-root}, it provides a geometric interpretation of its value. By a \emph{line} we mean a curve of degree one in $\PP^N$, where degrees of stable maps to $X$ and $Y$ are taken with respect to $H$. 

\begin{proposition}[Lines and the vanishing eigenvalue]\label{prop:lines}
Let $X\subset Y$ be a smooth hyperplane section of a smooth projective fivefold with $-K_Y=2H$, where $H$ is the hyperplane class. Suppose that the Fano schemes of lines $F(Y)$ and $F(X)$ are smooth of dimensions $4$ and $2$, and let $a_Y=\langle[\pt]\rangle^{Y}_{0,1,1}$ be the number of lines of $Y$ through a general point. Then
\begin{equation}\label{eq:lines}
  \langle v,w\rangle^{X}_{0,2,1}=-a_Y\,(v,w)_X\qquad\text{for all }v\in V_X\text{ and }w\in H^4(X,\QQ).
\end{equation}
In particular, if $H\star$ acts on $V$ by a scalar, then this scalar is $-a_Yq$.
\end{proposition}
 
\begin{proof}
Since $H$ is very ample, a degree one stable map with two marked points is a line $\ell\subset Y$ with $p_1,p_2\in\ell$. When these points collide, they move into a new component or ``bubble'', attached to the line at a node, and contracted by the map to a single point of $Y$. In either case,
\[
  \cM_Y\ceq\Mbar_{0,2}(Y,1)=P_Y\times_{F(Y)}P_Y,\qquad
  \cM_X\ceq\Mbar_{0,2}(X,1)=P_X\times_{F(X)}P_X,
\]
where $P\to F$ denotes the universal line. The stable maps with a bubble form the diagonal $D$ of $\cM_Y$, which is a copy of $P_Y=\Mbar_{0,1}(Y,1)$. By assumption $\cM_Y$, $\cM_X$ and $D$ are smooth of dimensions $6$, $4$ and $5$, which are the expected dimensions since $-K_Y=2H$ and $-K_X=H$. 
In particular, the evaluation map $\ev_D :\,D\to Y$ has degree $a_Y$. 

Let $\ev=(\ev_1,\ev_2)$ denote the evaluation maps and put
\[
  \cI_Y\ceq\ev_*[\cM_Y]\in H^8(Y\times Y,\QQ),\qquad
  \cI_X\ceq\ev_*[\cM_X]\in H^8(X\times X,\QQ).
\]
We claim that
\begin{equation}\label{eq:incidence}
  (i\times i)^*\cI_Y=\cI_X+a_Y\,[\Delta_X]\qquad\text{in }H^8(X\times X,\QQ),
\end{equation}
where $\Delta_X$ is the diagonal.  

A line in $Y$ either lies in $X$ or meets $X$ in a single point. Hence
$$
W\ceq \ev^{-1}({X\times X})=\cM_X\cup D_X, \qquad D_X\ceq\ev_D^{-1}({X})\subset D,
$$
and both pieces have dimension $4=\dim\cM_Y-2$. Each occurs with multiplicity one. Indeed, if $s\in H^0(Y,\cO_Y(1))$ defines $X$, then $W$ is cut out by $s(p_1)=s(p_2)=0$. Where $p_1\neq p_2$, these equations say that $s$ vanishes identically on $\ell$, so they cut out $\cM_X$, which is smooth of codimension
two. 

Near a line $\ell\not\subset X$, write $\sigma(\ell')=\ell'\cap X$ for nearby lines $\ell'$. Then $p_1=\sigma(\ell')$ and $p_2=\sigma(\ell')$, which again cut out a smooth subvariety of codimension two. The remaining locus consists of the bubbles on lines contained in $X$ and has dimension three. Thus
$$
  (i\times i)^![\cM_Y]=[\cM_X]+[D_X],
$$
 by \cite[Proposition~7.1]{Fulton1984}. Pushing forward by \cite[Theorem~6.2(a)]{Fulton1984} gives \eqref{eq:incidence}: the map $D_X\to X\times X$ factors through the diagonal $\Delta_X\cong X$ and the induced map $\ev_{D_X}\colon D_X\to X$ satisfies
$$
(\ev_{D_X})_*[D_X]=a_Y[X],
$$
by compatibility of proper pushforward with the refined pullback along $i$.

Now let $v\in V_X$ and $w\in H^4(X,\QQ)$. Since $V_X$ is orthogonal to $i^*H^4(Y,\QQ)$, we have $i_*v=0$, and so
\[
  \int_{X\times X}(v\boxtimes w)\smile(i\times i)^*\cI_Y
  =\int_{Y\times Y}(i_*v\boxtimes i_*w)\smile\cI_Y=0.
\]
By \eqref{eq:incidence} this reads $\langle v,w\rangle^X_{0,2,1}+a_Y\,(v,w)_X=0$, which proves \eqref{eq:lines}: on vanishing classes, the contribution of the lines of $X$ cancels that of the bubbles.
 
Finally, the coefficient of $q^d$ in $H\star v$ lies in $H^{6-2d}(X,\QQ)$, and by the divisor equation the coefficient of $q$ is
\[
  \sum_j\langle v,\phi_j\rangle^X_{0,2,1}\,\phi^j=-a_Y\,v,
\]
where $\{\phi_j\}$ and $\{\phi^j\}$ are dual bases of $H^4(X,\QQ)$. So if $H\star v=\lambda v$, then $\lambda=-a_Yq$.
\end{proof}
Gromov--Witten invariants and the decomposition $B_X\oplus V_X$ are invariant under deformations of the pair $(Y,X)$. So the conclusions of Proposition~\ref{prop:lines} hold as soon as $(Y,X)$ can be deformed to a pair for which $F(Y)$ and $F(X)$ are smooth of the expected dimensions.

\begin{corollary}[The cylinder correspondence]\label{cor:cylinder}
In the situation of Proposition~\ref{prop:lines}, let $p\,\colon P\to F$ be the universal line over $F=F(X)$ and $e\,\colon P\to X$ its evaluation map, and let $$
\Phi=p_*e^*\colon H^4(X,\QQ)\to H^2(F,\QQ)
$$ be the Abel--Jacobi map.
\begin{equation}\label{eq:cylinder}
  (\Phi v,\Phi w)_F=-a_Y\,(v,w)_X\qquad\text{for all } v\in V_X \text{ and }w\in H^4(X,\QQ).
\end{equation}
In particular, if $a_Y\neq0$, then $\Phi$ restricts to embeddings
$$\Phi: V_X(1)\hookrightarrow H^2(F,\QQ), \qquad \Phi: T_{X,\QQ}(1)\hookrightarrow T_{F,\QQ},$$
the second provided that $B_X \cap H^4(X,\QQ)$ consists of algebraic classes, so that $T_{X,\QQ}\subseteq V_X$.
\end{corollary}
 
\begin{proof}
Since $\Mbar_{0,2}(X,1)=P\times_FP$, the projection formula gives $\langle v,w\rangle^X_{0,2,1}=(\Phi v,\Phi w)_F$ for all $v,w\in H^4(X,\QQ)$, so \eqref{eq:cylinder} is \eqref{eq:lines}. It follows that $\Phi$ is injective on $V_X$ when $a_Y\neq0$, and $\Phi$ is a morphism of Hodge structures of bidegree $(-1,-1)$ because it is induced by an algebraic correspondence. Finally, if $v\in T_{X,\QQ}$ and $D\in\NS(F)\otimes\QQ$, then $(\Phi v,D)_F=(v,e_*p^*D)_X=0$, because $e_*p^*D$ is an algebraic class: so $\Phi$ maps $T_{X,\QQ}$ to $T_{F,\QQ}$.
\end{proof}
 
The map $\Phi$ is the transpose of the cylinder map $e_*p^*\colon H^2(F,\QQ)\to H^4(X,\QQ)$ studied by Lewis in~\cite{Lewis1985TheCH,Lew88} for hypersurfaces of degree $n$ in $\PP^n$; see also \cite{shimada}. For these hypersurfaces, where $a_Y=n!$, the eigenvalue $\lambda=-n!\,q$ is due to Collino-Jinzenji~\cite[Proposition~2]{Collino_1999}.

\section{Quantum Cohomology of K\"uchle fourfolds}\label{sec:comp}
In this section, we compute the small quantum Euler operator $U_X(0)=H\star$ of the K\"uchle fourfolds of types \kt{c1}, \kt{c6}, \kt{b8}, and \kt{b11}. In each case, Quantum Lefschetz gives the restriction of $U_X(0)$ to the inherited cohomology $B$, and the results of Section~\ref{sec:monodromy} then determine its action on the vanishing cohomology $V$. We briefly recall Quantum Lefschetz in \S\ref{subsec:ql} and outline the method in \S\ref{subsec:method}. 

\subsection{Quantum Lefschetz} \label{subsec:ql}
We use Quantum Lefschetz in the form of \cite[Theorem~1.1]{coates2014quantumlefschetzprinciplevector}, which relates Givental's Lagrangian cone of a zero locus to the Euler-twisted cone of the ambient variety. In this subsection only, the inclusion $i: X \hookrightarrow Y$ denotes the zero locus of a regular section of a convex vector bundle $\mathcal{E}$ on a smooth projective variety $Y$. 

The genus-zero Gromov--Witten invariants of $X$ are encoded by a Lagrangian cone $\mathcal{L}_X$ in the symplectic space 
$$
\mathcal H_X
=
H^\bullet(X;\Lambda_X^\text{Nov})\otimes\mathbb C((z^{-1})), \qquad 
\Omega_X(f,g)
=
\operatorname{Res}_{z=0}(f(-z),g(z))\,dz.
$$
Here $\Lambda_X^\text{Nov}$ is the \emph{Novikov ring} whose monomials $Q^\beta$ record degrees of effective curve classes $\beta\in \NE_\NN[X]\subset H_2(X,\ZZ)$. 

A general point of $L_X\subset \mathcal{H}_X$ has the form
\begin{equation}\label{eq:jfunc}
\mathbf J_X(\mathbf t) \ceq -z\mathbf1_X+\mathbf t(z) +\sum_{\beta, n, m, \mathbf{k}, \epsilon}\frac{Q^\beta}{n!}
\left\langle
t_{k_1}\psi_1^{k_1},\ldots,t_{k_n}\psi_n^{k_n},
\phi_\epsilon\psi_{n+1}^{m}
\right\rangle^X_{0,n+1,\beta}
\phi^\epsilon(-z)^{-m-1},
\end{equation}
where $\mathbf{t}(z)=t_0+t_1z+t_2z^2+\cdots$ with $t_k\in H^\bullet(X;\Lambda_X^{\mathrm{Nov}})$, and 
$\{\phi_\epsilon\}$ and $\{\phi^\epsilon\}$ are Poincar\'e dual bases of $H^\bullet(X;\mathbb{Q})$. The $J$-function $J_X(\tau, z)$ is the unique point of $\mathcal{L}_X$ of the form $-z + \tau +O(z^{-1}).$

Replacing the invariants of $X$ by the invariants of $Y$ twisted by the $S^1$-equivariant Euler class $\mathbf e$ of $\mathcal{E}$, with equivariant parameter $\kappa$, gives a cone $\mathcal L_{Y,\mathbf e}$ with a general point $\mathbf J_{Y,\mathbf e}$. It lives in the corresponding equivariant symplectic space for $Y$ with the twisted pairing 
$$
(a,b)_{Y,\mathbf e}
=\int_Y a\smile b\smile\mathbf e(\mathcal{E}).
$$

\begin{theorem}[Quantum Lefschetz, {\cite[Theorem~1.1]{coates2014quantumlefschetzprinciplevector}}]\label{thm:coates}
Let $i:X\hookrightarrow Y$ be the smooth zero locus of a regular section of a convex vector bundle $\mathcal{E}\to Y$. The non-equivariant limit of $\mathbf J_{Y,\mathbf e}(\mathbf t)$ exists and satisfies
\begin{equation}\label{eq:quantum_lefschetz}
\left.i^*\mathbf J_{Y,\mathbf e}(\mathbf t)\right|_{\kappa=0}
=\mathbf J_X(i^*\mathbf t).
\end{equation}
In particular,
$\left.i^*\mathcal L_{Y,\mathbf e}\right|_{\kappa=0} \subseteq\mathcal L_X$.
Throughout the homomorphism $Q^\delta \mapsto Q^{i{*}\delta}$ is applied to $\Lambda^{\text{Nov}}_X$.
\end{theorem}

Let $\{\phi_\alpha\}$ be a homogeneous basis of $H^\bullet(X, \QQ)$, with $\phi_\alpha \in H^{2d_{\alpha}}(X,\QQ)$, and write $\tau = \sum_\alpha \tau^\alpha \phi_\alpha$. We recover the big quantum product $\star_\tau$ from the $J$-function by Coates--Givental~\cite[Equation~16]{Coates2001QuantumRL}:
\begin{equation}\label{eq:16-cg}
 \begin{aligned}
 -z\frac{\partial^2 \mathbf{J}_X}{\partial\tau^\alpha\partial\tau^\beta}
 &=
\sum_\gamma A^\gamma_{\alpha,\beta}(\tau)
  \frac{\partial \mathbf{J}_X}{\partial\tau^\gamma},
  \qquad
  \phi_\alpha\star_\tau\phi_\beta
  &=
 \sum_\gamma A^\gamma_{\alpha,\beta}(\tau)\phi_\gamma.
  \end{aligned}
  \end{equation}
  
This formulation is particularly useful when $X$ inherits classes from $Y$ that are not powers of the hyperplane class, as in our examples.  It provides the input we need to compute the Hodge atom birational invariants of Katzarkov--Kontsevich--Pantev--Yu in~\cite{katzarkov2026birationalinvariantshodgestructures}. 

As discussed in the introduction, to use these invariants we study the generalized eigenspaces of quantum multiplication by the Euler vector field
$$
\operatorname{Eu}_X(\tau)
=c_1(X)+\sum_\alpha(1-d_\alpha)\tau^\alpha\phi_\alpha,
\qquad
U_X(\tau)(a)=\operatorname{Eu}_X(\tau)\star_\tau a.
$$
At $\tau=0$, the product $\star_0$ is the small quantum product of Section~\ref{sec:monodromy}. For the index-one, Picard-rank one fourfolds considered in this paper, we write $Q^\beta = q^{H\cdot \beta}$ and obtain
\begin{equation}\label{eq:small-euler}
U_X(0) = c_1(X)\star = H \star.
\end{equation}

\subsection{Method}\label{subsec:method}
Let $G=\Gr(k,n)$ be the Grassmannian of $k$-dimensional subspaces of $\CC^n$, with tautological subbundle $\mathcal S$, quotient bundle $\mathcal Q$, and Pl\"ucker line bundle $\mathcal O_G(1)=\det(\mathcal S^\vee)$. We write $Z_M(\mathcal F)$ for the zero locus of a general regular section of the vector bundle $\mathcal F$ on a variety $M$. Each fourfold considered below has the form
\begin{equation}
    X = Z_Y(\mathcal E), \qquad Y = Z_G(\mathcal F),
\end{equation}
where $\mathcal F$ is a globally generated homogeneous vector bundle on $G$ and $\mathcal E$ is a direct sum of copies of $\mathcal O_Y(1)$; see Table~2. 

By K\"uchle~\cite[Theorem~3]{kuchle}, $X$ is a Fano fourfold of index one, so that $$-K_X = H = c_1(\mathcal O_G(1))|_X = c_1(X).$$ 
It satisfies the standing assumptions of Section~\ref{sec:monodromy}, namely, it is Picard-rank one with $H^{\odd}(X,\QQ) = 0$. 

\begin{table}[H]\label{tab:kuchle-data}
    \centering
    \renewcommand{\arraystretch}{1.2}
    \setlength{\tabcolsep}{6pt}

    \begin{tabular}{@{}cclcccc@{}}
        \toprule
        K\"uchle Type
        & $G$
        & $\mathcal F$
        & $\mathcal E$
        & $h^{3,1}(X)$
        & $h^{2,2}(X)$
        & $b_4(X)$ \\
        \midrule
        \kt{b8}
        & $\Gr(2,7)$
        & $\operatorname{Sym}^2\mathcal S^\vee$
        & $\mathcal O_Y(1)^{\oplus3}$
        & $3$ & $38$ & $44$ \\

        \kt{b11}
        & $\Gr(2,8)$
        & $\mathcal Q^\vee(1)$
        & $\mathcal O_Y(1)^{\oplus2}$
        & $2$ & $31$ & $35$ \\
        
        \kt{c1}
        & $\Gr(3,6)$
        & ---
        & $\mathcal O_Y(1)^{\oplus5}$
        & $5$ & $52$ & $62$ \\

       \kt{c6}
        & $\Gr(3,7)$
        & $(\wedge^2\mathcal S^\vee)^{\oplus2}$
        & $\mathcal O_Y(1)^{\oplus2}$
        & $2$ & $30$ & $34$ \\

        \bottomrule
    \end{tabular}
    \medskip
    \caption{K\"uchle fourfolds from~\cite[Theorem~3]{kuchle} and their middle Hodge numbers. 
    }
\end{table}

Now let $Y'$ be a general five-dimensional linear section of $Y$ containing $X$. It is a smooth fivefold with hyperplane section $X$, and we apply Section~\ref{sec:monodromy} to $X\subset Y'$. By the Lefschetz hyperplane theorem $H^j(Y,\QQ)\to H^j(Y',\QQ)$ is an isomorphism
for $j\leq 4$, while $H^6(X,\QQ)$ and $H^8(X,\QQ)$ are spanned by $H^3$ and $H^4$.
Hence the inherited cohomology $B_X$ of Section~\ref{sec:monodromy} is the image of $H^\bullet(Y,\QQ)$.
We put $B_X^4\ceq B_X\cap H^4(X,\QQ)$, so that $V_X$ is the orthogonal complement of $B_X^4$ in $H^4(X,\QQ)$.

\begin{lemma}\label{lem:nondegen}
Suppose that $B_X^4$ is spanned by algebraic classes. Then
$$
h^{3,1}(V_X)=h^{3,1}(X),
\qquad
h^{2,2}(V_X)=h^{2,2}(X)-\dim B_X^4,
$$
and the Poincar\'e pairing on $V_X$ has signature
$\bigl(h^{2,2}(V_X),2h^{3,1}(V_X)\bigr)$.
\end{lemma}

\begin{proof}
Since $X$ is Fano, $h^{4,0}(X)=0$. The subspace $V_X$ is the orthogonal complement of $B_X^4$, which is of type $(2,2)$; hence it is a sub-Hodge structure of $H^4(X,\QQ)$ with the stated Hodge numbers.  

Let $v\in V_X$.  Since $H^2\in B_X^4$ we have $(H\smile v,H)=(v,H^2)=0$, and therefore $H\smile v=0$.  Thus $V_X$ consists of primitive classes, and Hodge--Riemann gives the stated signature relations.
\end{proof}

\begin{corollary}\label{cor:dimension}
Suppose that the middle inherited cohomology $B_X^4$ is spanned by algebraic classes and that $h^{3,1}(V_X)>0$ and $h^{2,2}(V_X)>0$.

Then $B$ is stable under $H\star$, and there is $\lambda\in K$ with the following properties.
\begin{enumerate}
\item
$H\star v=\lambda v$ for all $v\in V$, and $\lambda$ is a
repeated root of $\det(t-U_X(0)|_B)$.

\item
There is a unique class $b_\chi\in B$ such that $v\star w=(v,w)\,b_\chi$ for all $v,w\in V$.
It is isotropic and satisfies $H\star b_\chi=\lambda b_\chi$
and $(b_\chi,1)=1$.

\item
$\det(t-U_X(0))
=(t-\lambda)^{\dim V_X}\det(t-U_X(0)|_B)$.
\end{enumerate}
\end{corollary}

\begin{proof}
By Lemma~\ref{lem:nondegen}, the Poincar\'e pairing on $V_X$ is indefinite. Hence Proposition~\ref{prop:orthogonal-monodromy} applies to $X\subset Y'$, and the small quantum product is $O(V)$-equivariant. 
Proposition~\ref{prop:character} with $\lambda=\chi(H)$ shows that $B$ is invariant under $H\star$, and also 
gives (1) and (2): because the pairing on $V$ is nonzero, the class $b_\chi$ is unique. Part (3) follows from (1). 
\end{proof}

\textit{Inherited matrix computation.}
Assume that Corollary~\ref{cor:dimension} applies, and let $\mathcal B_X=(\varphi_1,\ldots,\varphi_s)$ be a homogeneous
basis of $B_X$, with $\varphi_j\in H^{2d_j}(X,\QQ)$. 

As $X$ has index one, $q$ has degree one, and the degree axiom gives
\begin{equation}\label{eq:daxiom}
H\star\varphi_j
=\sum_{i=1}^{s}a_{ij}\,q^{1+d_j-d_i}\varphi_i,
\qquad a_{ij}\in\QQ,
\end{equation}
with $a_{ij}=0$ whenever $d_i>d_j+1$.
We call $$A_X(q):=\bigl(a_{ij}q^{1+d_j-d_i}\bigr)_{i,j}$$ the \emph{inherited matrix} of $X$ in the basis $\mathcal B_X$; it is the matrix of $U_X(0)|_B$. 

Writing $D=\operatorname{diag}(q^{d_1},\ldots,q^{d_s})$, we have $A_X(q)=qD^{-1}A_X(1)D$, so that
\begin{equation}
\det\bigl(t-U_X(0)|_B\bigr)
=q^s\det\bigl((t/q)I_s-A_X(1)\bigr).
\end{equation}
Note that the entries of this matrix, for $d_i=d_j+1$, are the structure constants of the cup product with $H$. Moreover, $H\star$ is self-adjoint for the Poincar\'e pairing, so $G_XA_X(q)$ is symmetric, where $G_X$ is the Gram matrix of $\mathcal B_X$.

To compute $A_X(q)$, we apply Theorem~\ref{thm:coates} to a family $I_{Y,e}$ of points of the twisted cone from Section~\ref{subsec:ql}. 
This gives
$$
I_X(q,z)
=i^*I_{Y,e}(q,z;\kappa)|_{\kappa=0}
\in\mathcal L_X\subset \mathcal{H}_X,
$$
and in each case $I_X(q,z)=-z+a_Xq+O(z^{-1})$ for some $a_X\in\QQ$.
By the string equation and the characterization of $J_X$ in \S\ref{subsec:ql},
\begin{equation}\label{eq:string-normalization}
I_X(q,z)
=J_X(a_Xq\,\mathbf 1;q,z)
=e^{-a_Xq/z}J_X(0;q,z).
\end{equation}
Writing $\tau=\sum_i t^i\varphi_i$ and applying \eqref{eq:16-cg} to our normalization of the family $J_X(\tau;q,z)$ at $\tau=0$ gives $A_X(q)$.
Equivalently, if $\widetilde A_X(q)$ is the matrix of $U_X(a_Xq\,\mathbf 1)|_B$, then
\begin{equation}\label{eq:string-matrix-shift}
A_X(q)=\widetilde A_X(q)-a_XqI_s.
\end{equation}

\begin{remark}[String shift and the vanishing eigenvalue]\label{rem:shift}
Let $a_X$ be the coefficient of the unit-coordinate correction in \eqref{eq:string-normalization}. We show that it is the number of lines of $Y'$ through a general point. Indeed, the hypergeometric modification for the last hyperplane section multiplies the degree-$d$ term of $J_{Y'}(0;q,z)$ by
$\prod_{m=1}^{d}(H-mz)$. Since $Y'$ has index two, the coefficient of $q$ in $J_{Y'}(0;q,z)$ is $-a\,z^{-1}\bone+O(z^{-2})$ with $a=\langle[\pt]\rangle^{Y'}_{0,1,1}$, and multiplying by $H-z$ gives
\begin{equation}\label{eq:linecount}
 a_X=\langle[\pt]\rangle^{Y'}_{0,1,1}. 
\end{equation}
Thus the unit-coordinate correction is $a_Xq\bone$, which agrees with Iritani's mirror-map formula for a hyperplane section of an index-two Fano variety \cite[Section~0.2]{Iritani2017QuantumLefschetz}.

Because \eqref{eq:linecount} is deformation invariant, this holds for every smooth pair $X\subset Y'$ of the given type. For types \kt{b8}, \kt{b11}, \kt{c1}, \kt{c6} we have $a_X=4,3,6,5$\footnote{The method for deriving the stated line counts is provided in
Appendix~\ref{app:lines}.}, which agree with the vanishing eigenvalues $\lambda=-a_Xq$ computed in \S\ref{sec:comp}. 

Consequently, the translation in \eqref{eq:string-matrix-shift} gives an operator annihilating the vanishing summand:
$$
\left.U_X(a_Xq\,\mathbf 1)\right|_{V_X\otimes K}=0\qquad \text{so that}\qquad \det\bigl(tI-U_X(a_Xq\,\mathbf 1)\bigr)
=
t^{\dim V_X}\det\bigl(tI_s-\widetilde A_X(q)\bigr).
$$
This gives a cohomological interpretation of what Golyshev refers to as the \emph{Givental constant} in \cite[Section~1.2]{golyshev2008spectrastrains}. There, he explains that the Quantum Lefschetz comparison of spectra requires an additional translation when the hyperplane section has index one. After this translation, the vanishing cohomology contributes only zero eigenvalues to the shifted operator, which is the behavior suggested by Iritani's question on extending the Quantum Lefschetz comparison of characteristic polynomials to include primitive cohomology \cite[Question~15]{Iritani2017QuantumLefschetz}. 
\end{remark}

\subsection{Applications to K\"uchle fourfolds}\label{subsec:kuchlecomp}
\subsubsection{K\"uchle fourfold of type \kt{b8}}\label{subsec:b8}
Let $G=\Gr(2,7)$ and 
$$
Y=Z_G\!\left(\operatorname{Sym}^2\mathcal S^\vee\right),
\qquad
\mathcal E=\mathcal O_Y(1)^{\oplus3},
\qquad
X=Z_Y(\mathcal E).
$$

Put $P=c_2(\mathcal S^\vee)|_X$. The inherited cohomology $B_X$ has the basis
$$
\mathcal B_{b8}=(1,H,H^2,P,H^3,H^4),
$$
and its Poincar\'e pairing is determined by the intersection numbers 
$$
H^4 = 56, \qquad H^2\cdot P=20, \qquad P^2 = 8. 
$$
Since the classes $H^2$ and $P$ are algebraic, Lemma~\ref{lem:nondegen} gives $h^{3,1}(V_X) =3$ and $h^{2,2}(V_X) =36$. In particular, $\dim V_X = 42$ and Corollary~\ref{cor:dimension} applies. 

\begin{proposition}[Type \kt{b8}]\label{prop:quantum-b8}
Let $X$ be a K\"uchle fourfold of type \kt{b8}. 
\begin{enumerate}
    \item  In the basis $\mathcal B_{b8}$, the matrix of $U_X(0)|_B$ is
$$
A_{b8}(q)=
\begin{pmatrix}
0 & 96q^2 & 2592q^3 & 936q^3 & 51968q^4 & 434560q^5\\
1 & 6q & 300q^2 & 108q^2 & 6048q^3 & 51968q^4\\
0 & 1 & 6q & 4q & 280q^2 & 2352q^3\\
0 & 0 & 8q & -2q & 56q^2 & 672q^3\\
0 & 0 & 1 & \tfrac{5}{14} & 6q & 96q^2\\
0 & 0 & 0 & 0 & 1 & 0
\end{pmatrix}$$
Consequently, 
$$ 
\det(t - U_X(0)|_B) =(t+4q)^3(t+8q)(t^2-36qt-176q^2).
$$    
\item For all $v,w\in V$,
$$
H\star v = -4qv, \qquad\text{and}\qquad v\star w = (v,w) b_\chi,
$$
where $b_\chi=b/56$ and the class
$$
b=H^4-4qH^3-56q^2H^2-672q^3H-4256q^4
$$
spans the $(-4q)$-eigenspace of $U_X(0)|_B$. 
\item We have 
$$
\det(t - U_X(0))=(t+4q)^{45}(t+8q)(t^2-36qt-176q^2). 
$$
In particular $M_X = 45$.
\end{enumerate}
\end{proposition}

\begin{proof}
(1) Apply the method of \S\ref{subsec:method} to the Euler-twisted Grassmannian $I$-function of $G$, twisted by $\operatorname{Sym}^2 S^\vee\oplus\mathcal O_G(1)^{\oplus 3}$, constructed using the abelian/non-abelian correspondence of Inoue--Ito--Miura~\cite[Definition~3.3, Theorem~3.4]{Inoue_2017}. This gives $A_{b8}(q)$ and the characteristic polynomial follows. 

(2) By Proposition~\ref{prop:character}, $H\star$ acts on $V$ by a scalar which is a repeated root of the polynomial from the previous step. Hence $\lambda=-4q$. A direct computation with $A_{b8}(q)$ shows that the $(-4q)$-eigenspace of $U_X(0)|_B$ is the line spanned by $b$
and $(b,1)=H^4=56$. By Proposition~\ref{prop:character}, this line contains $b_\chi$, and $(b_\chi,1)=1$. Hence $b_\chi=b/56$.
  
(3) This is Corollary~\ref{cor:dimension} since $\dim V_X = 42$.
\end{proof}

\subsubsection{K\"uchle fourfold of type \kt{b11}}
Let $G=\Gr(2,8)$ and 
$$
\qquad
Y=Z_G\!\left(\mathcal Q^\vee(1)\right),
\qquad
\mathcal E=\mathcal O_Y(1)^{\oplus2},
\qquad
X=Z_Y(\mathcal E).
$$

Let $\mathcal U_1,\mathcal U_2$ be Meng's universal bundles on $Y$, of ranks two and three~\cite{meng2024quantumcohomologyfanoquiver}, and put
$$
c_2=c_2(\mathcal U_2^\vee),\qquad
d_2=c_2(\mathcal U_1^\vee),\qquad
H_Y=c_1(\mathcal U_1^\vee)=c_1(\mathcal U_2^\vee).
$$
Define
$$
e_1=-c_2+2d_2,\qquad
e_2=H_Y^2-2d_2,\qquad
e_3=c_2-d_2,\qquad
\bar e_j=e_j|_X.
$$
The inherited cohomology $B_X$ has the basis
$$
\mathcal B_{b11}=(1,H,\bar e_1,\bar e_2,\bar e_3,H^3,\text{pt}),\qquad  \text{pt} \ceq \frac{H^4}{57}
$$
and its Poincar\'e pairing is determined by $H^4=57$, the relation $H^2 = 2 \bar e_1 + \bar e_2 + 2 \bar e_3$ and the Gram matrix of $(\bar e_1, \bar e_2, \bar e_3)$, which is  
$$
\begin{pmatrix}
2&3&1\\
3&9&3\\
1&3&2
\end{pmatrix}.
$$ 
Since the classes $\bar e_j$ are algebraic, Lemma~\ref{lem:nondegen} gives $h^{3,1}(V_X) =2$ and $h^{2,2}(V_X) =28$. In particular, $\dim V_X = 32$ and Corollary~\ref{cor:dimension} applies. 
\begin{proposition}[Type \kt{b11}]\label{prop:quantum-b11}
Let $X$ be a K\"uchle fourfold of type \kt{b11}. 
\begin{enumerate}
    \item  In the basis $\mathcal B_{b11}$, the matrix of $U_X(0)|_B$ is
$$
A_{b11}(q)=\begin{pmatrix}
0&84q^2&333q^3&846q^3&333q^3&42408q^4&5805q^5\\
1&5q&42q^2&108q^2&42q^2&5301q^3&744q^4\\
0&2&-2q&3q&2q&342q^2&51q^3\\
0&1&3q&6q&3q&456q^2&60q^3\\
0&2&2q&3q&-2q&342q^2&51q^3\\
0&0&3/19&7/19&3/19&5q&28q^2/19\\
0&0&0&0&0&57&0
\end{pmatrix}.
$$
Consequently, 
$$ 
\det(t - U_X(0)|_B) =(t+3q)^3(t+4q)^2(t^2-29qt-339q^2).
$$    
\item For all $v,w\in V$,
$$
H\star v = -3qv, \qquad\text{and}\qquad v\star w = (v,w) b_\chi,
$$
where $b_\chi=b/57$ and the class
$$
b=H^4-3qH^3
  -57q^2(\bar e_1+2\bar e_2+\bar e_3)
  -684q^3H-3933q^4
$$
spans the $(-3q)$-eigenspace of $U_X(0)|_B$. 
\item We have 
$$
\det(t - U_X(0))=(t+3q)^{35}(t+4q)^2(t^2-29qt-339q^2). 
$$
In particular $M_X = 35$.
\end{enumerate}
\end{proposition}
\begin{proof}

 (1) We obtain $J_Y$ from Meng's quantum Chevalley formula \cite[Corollary~5.26]{meng2024quantumcohomologyfanoquiver}, and apply Coates--Givental hypergeometric modification \cite[Theorem~2]{Coates2001QuantumRL} for $\mathcal O_Y(1)^{\oplus2}$. The method of \S\ref{subsec:method} then gives $A_{b11}(q)$ and the characteristic polynomial follows.
 
(2) By Corollary~\ref{cor:dimension}, $H\star$ acts on $V$ by a scalar $\lambda$, which is a repeated root of the polynomial in (1). The repeated roots are $-3q$ and $-4q$: the factor $t^2-29qt-339q^2$ has discriminant $2197q^2$, and does not vanish at $t=-3q$ nor at $t=-4q$.  

We exclude $-4q$ because a computation with $A_{b11}(q)$ shows that the $(-4q)$-eigenspace of $U_X(0)|_B$ is the plane spanned by
$$
\xi=\bar e_1-\bar e_3\qquad\text{and}\qquad
\eta=2H^4-8qH^3-57q^2(3\bar e_1+2\bar e_2+3\bar e_3)
-1026q^3H-6498q^4,
$$
and that
$$
(\xi,\xi)=2,\qquad (\xi,\eta)=0,\qquad
(\eta,\eta)=-6\,(57q^2)^2.
$$
Hence, for $\alpha,\beta\in K$,
$$
(\alpha\xi+\beta\eta,\ \alpha\xi+\beta\eta)
=2\bigl(\alpha^2-3\,(57q^2\beta)^2\bigr),
$$
which vanishes only for $\alpha=\beta=0$ (because $3$ is not a square in $K$). Therefore the eigenspace contains no nonzero isotropic vector, whereas $b_\chi$ is isotropic. Hence $\lambda=-3q$.

A computation with $A_{b11}(q)$ shows that the $(-3q)$-eigenspace of $U_X(0)|_B$ is the line spanned by $b$,
and $(b,1)=H^4=57$. As this line contains $b_\chi$ and $(b_\chi,1)=1$, we get $b_\chi=b/57$.
    
(3) This is Corollary~\ref{cor:dimension} since $\dim V_X = 32$.
\end{proof}

\subsubsection{K\"uchle fourfold of type \kt{c1}}\label{subsec:c1}
Let $Y=G=\Gr(3,6)$ and 
$$
\qquad
\mathcal E=\mathcal O_Y(1)^{\oplus5},
\qquad
X=Z_Y(\mathcal E).
$$
Put $P=c_2(S^\vee)|_X$. The inherited cohomology $B_X$ has the basis
$$
\mathcal B_{c1}=(1,\ H,\ H^2,\ P,\ H^3,\ H^4),
$$
and its Poincar\'e pairing is determined by the intersection numbers
$$
H^4=42,\qquad H^2\cdot P=21,\qquad P^2=11.
$$
Since $H^2$ and $P$ are algebraic, Lemma~\ref{lem:nondegen} gives $h^{3,1}(V_X)=5$ and $h^{2,2}(V_X)=50$. In particular $\dim V_X=60$, and Corollary~\ref{cor:dimension} applies.

\begin{proposition}[Type \kt{c1}]
Let $X$ be a K\"uchle fourfold of type \kt{c1}.
\begin{enumerate}
\item In the basis $\mathcal B_{c1}$, the matrix of $U_X(0)|_B$ is
$$
A_{c1}(q)=
\begin{pmatrix}
0 & 180q^2 & 6804q^3 & 3402q^3 & 184464q^4 & 2026080q^5 \\
1 & 10q    & 588q^2  & 294q^2  & 16380q^3  & 184464q^4 \\
0 & 1      & 15q     & 10q     & 588q^2    & 6804q^3 \\
0 & 0      & 0       & -5q     & 0         & 0 \\
0 & 0      & 1       & \frac12 & 10q       & 180q^2 \\
0 & 0      & 0       & 0       & 1         & 0
\end{pmatrix}.
$$
Consequently
$$
\det(t-U_X(0)|_B)=(t+6q)^3(t+5q)^2(t-58q).
$$
\item For all $v,w\in V$,
$$
H\star v=-6q\,v
\qquad\text{and}\qquad
v\star w=(v,w)\,b_\chi,
$$
where $b_\chi = b/42$ and the class
$$
b=H^4-6qH^3-84q^2H^2-1512q^3H-12600q^4
$$
spans the $(-6q)$-eigenspace of $U_X(0)|_B$.
\item We have
$$
\det(t-U_X(0))=(t+6q)^{63}(t+5q)^2(t-58q).
$$
In particular $M_X=63$.
\end{enumerate}
\end{proposition}

\begin{proof}
(1) Apply the method of \S\ref{subsec:method} to the Euler-twisted $I$-function of $G$, twisted by $\mathcal O_G(1)^{\oplus 5}$, using the abelian/non-abelian correspondence of Inoue--Ito--Miura~\cite[Definition~3.3, Theorem~3.4]{Inoue_2017}. This gives $A_{c1}(q)$, and the
characteristic polynomial follows.

(2) By Corollary~\ref{cor:dimension}, $H\star$ acts on $V$ by a scalar $\lambda$ which is a repeated root of the polynomial in (1), so $\lambda\in\{-6q,-5q\}$. We exclude $-5q$ as in the proof of Proposition~\ref{prop:quantum-b11}: a computation with $A_{c1}(q)$ shows that the $(-5q)$-eigenspace of $U_X(0)|_B$ is the plane spanned by
$$
  \zeta=H^2-2P\qquad\text{and}\qquad\zeta'=H^4-5qH^3-105q^2H^2-1764q^3H-14364q^4,
$$
and that $(\zeta,\zeta)=2$, $(\zeta,\zeta')=0$ and $(\zeta',\zeta')=-6\,(21q^2)^2$. Since $3$ is not a square in $K$, this plane contains no nonzero isotropic vector, whereas $b_\chi$ is isotropic. Hence $\lambda=-6q$, in agreement with Remark~\ref{rem:shift}, where $a_X=6$. The $(-6q)$-eigenspace of $U_X(0)|_B$ is the line spanned by $b$
and $(b,1)=H^4=42$. As this line contains $b_\chi$ and $(b_\chi,1)=1$, we get $b_\chi=b/42$.

(3) This is Corollary~\ref{cor:dimension} since $\dim V_X=60$.
\end{proof}

\subsubsection{K\"uchle fourfold of type \kt{c6}}
Let $G=\Gr(3,7)$ and 
$$
\qquad
Y=Z_G\!\left(
\left(\bigwedge\nolimits^2\mathcal S^\vee\right)^{\oplus2}
\right),
\qquad
\mathcal E=\mathcal O_Y(1)^{\oplus2},
\qquad
X=Z_Y(\mathcal E).
$$

Put $P=c_2(\mathcal Q)|_X$. The inherited cohomology $B_X$ has the basis
$$
\mathcal B_{c6}=(1,\ H,\ H^2,\ P,\ H^3,\ H^4),
$$
and its Poincar\'e pairing is determined by the intersection numbers
$$
H^4=61,\qquad H^2\cdot P=33,\qquad P^2=18.
$$
Since $H^2$ and $P$ are algebraic, Lemma~\ref{lem:nondegen} gives $h^{3,1}(V_X)=2$ and $h^{2,2}(V_X)=28$. In particular $\dim V_X=32$, and Corollary~\ref{cor:dimension} applies.

\begin{proposition}[Type \kt{c6}]
Let $X$ be a K\"uchle fourfold of type \kt{c6}.
\begin{enumerate}
\item In the basis $\mathcal B_{c6}$, the matrix of $U_X(0)|_B$ is
$$
A_{c6}(q)=
\begin{pmatrix}
0 & 88q^2 & 2097q^3 & 1134q^3       & 38308q^4 & 311710q^5 \\
1 & 6q    & 244q^2  & 132q^2        & 4575q^3  & 38308q^4 \\
0 & 1     & 13q     & 9q            & 244q^2   & 2196q^3 \\
0 & 0     & -8q     & -8q           & 0        & -183q^3 \\
0 & 0     & 1       & \frac{33}{61} & 6q       & 88q^2 \\
0 & 0     & 0       & 0             & 1        & 0
\end{pmatrix}.
$$
Consequently
$$
\det(t-U_X(0)|_B)
=(t+5q)^3\bigl(t^3-32qt^2-195q^2t-293q^3\bigr).
$$
\item For all $v,w\in V$,
$$
H\star v=-5q\,v
\qquad\text{and}\qquad
v\star w=(v,w)\,b_\chi,
$$
where $b_\chi = b/61$ and the class
$$
b= H^4-5qH^3-61q^2P-427q^3H-2684q^4
$$
spans the $(-5q)$-eigenspace of $U_X(0)|_B$.
\item We have
$$
\det(t-U_X(0))
=(t+5q)^{35}\bigl(t^3-32qt^2-195q^2t-293q^3\bigr).
$$
In particular $M_X=35$.
\end{enumerate}
\end{proposition}
\begin{proof}
(1) Apply the method of \S\ref{subsec:method} to the Euler-twisted Grassmannian $I$-function of $G$, twisted by $\left(\bigwedge\nolimits^2\mathcal S^\vee\right)^{\oplus2}\oplus\mathcal O_G(1)^{\oplus2}$, using the  abelian/non-abelian correspondence of Inoue--Ito--Miura~\cite[Definition~3.3, Theorem~3.4]{Inoue_2017}. This gives $A_{c6}(q)$, and the
characteristic polynomial follows.
 
(2) By Corollary~\ref{cor:dimension}, $H\star$ acts on $V$ by a scalar $\lambda$, which is a repeated root of the polynomial in (1). Its cubic factor is squarefree over $\QQ(q)$, and gives the nonzero value $-243q^3$ at $t=-5q$. Thus $-5q$ is the unique repeated root and $\lambda=-5q$. 

A computation with $A_{c6}(q)$ shows that the $(-5q)$-eigenspace of $U_X(0)|_B$ is the line spanned by $b$ and $(b,1)=H^4=61$. As this line contains $b_\chi$ and $(b_\chi,1)=1$, we get $b_\chi=b/61$.

(3) This is Corollary~\ref{cor:dimension} since $\dim V_X=32$.
\end{proof}

The following statement is used in Section~\ref{sec:results}.
\begin{proposition}\label{prop:c6-saturation}
Let $X\subset Y \subset G$ be a K\"uchle fourfold of type \kt{c6} and $P= c_2(\mathcal Q)|X$. Then the classes $H^2$ and $P$ generate a saturated sublattice of $H^4(X,\ZZ)/\mathrm{tors}$.
\end{proposition}
\begin{proof}
Let $\omega_0,\omega_1\in \bigwedge^2(\CC^7)^\vee$ be the alternating $2$-forms defining the smooth sixfold $Y\subset G$, and let $E \subset \CC^7$ be a common $4$-dimensional isotropic subspace for $\omega_0$ and $\omega_1$ (any pair of alternating $2$-forms on $\CC^{2m+1}$ admits a common $(m+1)$-dimensional isotropic subspace). 

Choose $w \in \CC^7 \setminus E$, a line $\ell \subset E$, and a $2$-plane $A \subset E \cap w^{\perp_{\omega_0}} \cap w^{\perp_{\omega_1}}$. Then the Schubert planes
$$\Pi_0 = \{ L \in G \mid \ell \subset L \subset E \} \cong \PP^2,
  \qquad
  \Pi_1 = \{ L \in G \mid A \subset L \subset E \oplus \langle w \rangle \} \cong \PP^2$$

lie in $Y$. Writing $H_Y = c_1(\mathcal{O}_Y(1))$, standard Schubert calculus gives 
\[
  \begin{pmatrix}
    \int_{\Pi_0} H_Y^2 & \int_{\Pi_0} c_2(\mathcal{Q}) \\[3pt]
    \int_{\Pi_1} H_Y^2 & \int_{\Pi_1} c_2(\mathcal{Q})
  \end{pmatrix}
  =
  \begin{pmatrix}
    1 & 0 \\
    1 & 1
  \end{pmatrix}.
\]
Because this pairing matrix is unimodular, $H_Y^2$ and $c_2(\mathcal{Q})|_Y$ span a saturated rank-$2$ sublattice of $H^4(Y,\ZZ)/\mathrm{tors}$. Since $X$ is cut out inside $Y$ by two hyperplane sections, integral Weak Lefschetz  shows that $H^4(Y,\ZZ) \hookrightarrow H^4(X,\ZZ)$ is injective with torsion-free cokernel, so saturation passes to $H^4(X,\ZZ)/\mathrm{tors}$.
\end{proof}

\section{Quantum Noether Rationality criterion}\label{sec:lattices}
In \cite[Proposition~26]{Gueintegral}, Gu\'er\'e proves an integral Hodge isometry between the transcendental part of a Hodge atom of a rational fourfold and the transcendental lattice of an associated minimal surface. Here we combine it with the \emph{Hodge atoms} of~\cite{katzarkov2026birationalinvariantshodgestructures} to prove Proposition~\ref{pro:quantum-integral}. To do this, we consider the generalized eigenspaces of $U_X(\tau)$ at general points of the
non-archimedean Hodge locus of \cite[Section~5.4]{katzarkov2026birationalinvariantshodgestructures}.

We first recall the elementary lattice theory used below; see Nikulin~\cite[Section~1]{nikulin}.

\subsection{Lattices}
From now on, \emph{lattice} means a free abelian group $L$ of finite rank equipped with a symmetric bilinear pairing
$(\cdot, \cdot):L\times L\longrightarrow\ZZ$ that is nondegenerate over $\QQ$. Write
$$
L^\vee=\{x\in L\otimes\QQ:(x,L)\subset\ZZ\},\qquad
D(L)=L^\vee/L,
$$
and let $\ell(D(L))$ denote the minimum number of generators of the discriminant group $D(L)$. A lattice $L$ is unimodular if $L=L^\vee$; the lattices $\Lambda_X$ and $\Lambda_S$ are unimodular by integral Poincar\'e duality. A sublattice $L\subseteq\Lambda$ is \emph{primitive} (or \emph{saturated}) if $\Lambda/L$ is torsion-free, equivalently if $L=(L\otimes\QQ)\cap\Lambda$.

If $\Gamma_L$ is the Gram matrix of $L$ in an integral basis, then
$$
D(L)\cong\operatorname{coker}\Gamma_L,\qquad
|D(L)|=|\det\Gamma_L|,\qquad
\ell(D(L))\leq\rk L.
$$
Thus the discriminant group can be computed from the Smith normal form of $\Gamma_L$.

If $L\subset M$ are lattices of the same rank, then $M/L$ is a finite group. Its number of elements, $[M:L]\ceq|M/L|$, is called the \emph{index} of $L$ in $M$. Since their Gram matrices satisfy
$$
\det\Gamma_L=[M:L]^2\det\Gamma_M,
$$
a sublattice is saturated whenever $|\det\Gamma_L|$ is squarefree.

If $L$ is primitive and nondegenerate and $\Lambda$ is unimodular, then
$D(L)\cong D(L^\perp),$ as finite abelian groups~\cite[Proposition~1.6.1]{nikulin}.

\begin{proof}[Proof of Proposition~\ref{pro:quantum-integral}]
Suppose that $X$ is rational. By Gu\'er\'e's integral Hodge isometry \cite[Propositions~25--26]{Gueintegral}, applied to a weak factorization of a birational map $X\dashrightarrow\PP^4$, there is a surface centre whose smooth minimal model $S$ satisfies 
$$T_X(1)\simeq T_S$$ 
as integral Hodge structures. The isometry gives $p_g(S)=h^{3,1}(X)>0$. Thus $S$ is not ruled and minimality implies that $K_S$ is nef. 

Since $T_{X,\QQ}$ is polarizable and $\End_{\Hdg}(T_{X,\QQ})=\QQ\cdot \id$, it is a simple Hodge structure and it is not of Tate type (as $h^{3,1}(X)>0$), where its orthogonal complement in $H^\even(X,\QQ)$ is entirely Tate by \eqref{eq:conditions}. Given that Gromov--Witten classes are algebraic, it follows that $U_X(\tau)$ preserves $T_X \otimes K$ and acts on it by a scalar $\lambda_X(\tau)$. In particular, $T_{X,\QQ}$ lives in a single Hodge atom of $X$, which corresponds (under weak factorization) to an atom of $S$ and all of $H^\bullet(S,\QQ)$; see \cite[Lemma~5.24]{katzarkov2026birationalinvariantshodgestructures}. After specialization this gives
\begin{equation}
\label{eq:atom-bound}
r+\rho(S)+2+4h^{1,0}(S)
=\dim_{\QQ}H^\bullet(S,\QQ)
\leq m_X\leq M_X.
\end{equation}

The lattices $T_S$ and $\NS(S)\tors$ are primitive
orthogonal complements in the unimodular lattice $\Lambda_S$.
Hence 
$$
D(T_X)\simeq D(T_S)
\simeq D(\NS(S)\tors),
\qquad
\ell(D(T_X))\leq\rho(S).
$$
Together with $\rho(S)\geq1$ and \eqref{eq:atom-bound}, this proves \textup{(1)} and the lower bound in \textup{(2)}. Noether's formula gives
\begin{equation}
\label{eq:noether-S}
K_S^2+\rho(S)+8h^{1,0}(S)
=10+12h^{3,1}(X)-r.
\end{equation}
Since $K_S^2\geq0$, the upper bound in \textup{(2)} follows.

Finally, suppose $M_X\leq r+3$.
Equation~\eqref{eq:atom-bound} requires
$$
\rho(S)=1,\qquad h^{1,0}(S)=0,\qquad
m_X=M_X=r+3.
$$
Consequently, $e(S)=M_X$ and
$$
\NS(S)\tors=\ZZ u,
\qquad u^2=|D(T_X)|,
$$
so $D(T_X)$ is cyclic.
Writing $[K_S]=ku$ modulo torsion, with $k\in\ZZ$, Noether's formula becomes
$$
12\bigl(1+h^{3,1}(X)\bigr)
=12\chi(\mathcal O_S)
=K_S^2+e(S)
=|D(T_X)|k^2+M_X.
$$
This proves \textup{(3)}.
\end{proof}

\section{Irrationality of K\"uchle fourfolds}\label{sec:results}
\begin{proposition}\label{prop:hodge-lattices}
For a very general K\"uchle fourfold $X$ of type \kt{b8}, \kt{b11}, \kt{c1} or \kt{c6},
$$
T_{X,\QQ}=V_X
\qquad\text{and}\qquad
\End_{\mathrm{Hdg}}(T_{X,\QQ})=\QQ\cdot\operatorname{id}.
$$
Moreover, with $P$ as in Section~\ref{sec:comp}, a $\ZZ$-basis of $A_X$ is given by 
$\bigl(P,\tfrac12c_2(\mathcal Q)|_X\bigr)$ in type \kt{b8}, 
$(\bar e_1,\bar e_2,\bar e_3)$ in type \kt{b11},  and by
$(H^2,P)$ in types \kt{c1} and \kt{c6}, with Gram matrices
$$
\Gamma_{c1}=
\begin{pmatrix}
42&21\\
21&11
\end{pmatrix},
\qquad
\Gamma_{c6}=
\begin{pmatrix}
61&33\\
33&18
\end{pmatrix},
\qquad
\Gamma_{b8}=
\begin{pmatrix}
8&6\\
6&6
\end{pmatrix},
\qquad
\Gamma_{b11}=
\begin{pmatrix}
2&3&1\\
3&9&3\\
1&3&2
\end{pmatrix}.
$$
In particular $D(T_X)\simeq D(A_X)$ is isomorphic to $\ZZ/2\oplus\ZZ/6$ and $(\ZZ/3)^2$, $\ZZ/21$, $\ZZ/9$, respectively.
\end{proposition}
\begin{proof}
For very general $X$, the Hodge group of $V_X$ contains the connected algebraic monodromy group $SO(V_X)$; see \cite[Proposition~7.5]{Deligne1972}.
Its standard representation is absolutely irreducible and contains no nonzero invariant vectors. Thus $V_X$ contains no nonzero rational Hodge classes and Schur's lemma gives:
$$
\operatorname{End}_{\mathrm{Hdg}}(V_X)
=\QQ\cdot\operatorname{id}.
$$
Since $B_X^4$ is algebraic, it follows that
$$
A_X\otimes_{\ZZ}\QQ=B_X^4,
\qquad
T_{X,\QQ}=V_X.
$$

Let $\mathcal{L}_X$ be the lattice generated by these classes, whose rational span is $B_X^4$ by Section~\ref{subsec:kuchlecomp}. 
\begin{itemize}
    \item For \kt{c1}, $\det\Gamma_{c1}=21$ is squarefree, so $\mathcal{L}_X$ is saturated.
    \item For \kt{b8}, standard integral Schubert calculus on $Y=\operatorname{OG}(2,7)$ gives
$$
H^4(Y,\ZZ)=\ZZ H_Y^2\oplus\ZZ\tau_2,
\qquad c_2(\mathcal Q)|_Y=2\tau_2.
$$ Thus $(c_2(\mathcal S^\vee)|_Y,\tau_2)$ is also an integral basis. 
\item For \kt{b11}, the classes $e_1,e_2,e_3$ form an integral basis by Meng's cellular calculation
\cite[Sections~3.3--4, Proposition~4.4]{meng2024quantumcohomologyfanoquiver}.
\item For \kt{c6}, saturation follows from Proposition~\ref{prop:c6-saturation}. 
\end{itemize}

In types \kt{b8} and \kt{b11}, integral weak Lefschetz through the successive smooth hyperplane sections shows that $\mathcal{L}_X$ is primitive in $\Lambda_X$ as well. Since its rational span is $A_X\otimes_{\ZZ}\QQ$, this gives $L_X= A_X$.

The displayed Gram matrices are the same as those from \S\ref{sec:comp}, except for \kt{b8}, where $H^2=P+2\theta$ with $\theta=\tfrac12c_2(\mathcal Q)|_X$ and we use the basis $(P,\theta)$. The stated discriminant groups follow from their Smith normal forms, and unimodularity of $\Lambda_X$ identifies $D(T_X)$ with $D(A_X)$.

\end{proof}

\begin{proof}[Proof of Theorem~\ref{cor:kuchleirr}]
Let $X$ be very general of one of the four types. By Proposition~\ref{prop:hodge-lattices}, $X$ satisfies \eqref{eq:conditions} with $T_{X,\QQ}=V_X$. 
The computations in \S\ref{subsec:kuchlecomp} show that we have $M_X=\rk T_X+3$ in each case.

Suppose that $X$ is rational.  Proposition~\ref{pro:quantum-integral}(3) then requires $D(T_X)$ to be cyclic. Hence types \kt{b8} and \kt{b11} cannot be rational by Proposition~\ref{prop:hodge-lattices}. 
Moreover, $|D(T_X)|\,k^2=9+12h^{3,1}(X)-\rk T_X$ for some $k\in\ZZ$. For types \kt{c1} and \kt{c6} this is:
$$
    21k^2=9, \qquad 9k^2=1,  
$$
respectively, and neither has an integer solution. 

\end{proof}

\section{Applications to Verra fourfolds}\label{sec:verra}
We now apply Proposition~\ref{pro:quantum-integral} to recover Gu\'er\'e's result that the very general Verra fourfold is irrational~\cite[Corollary~28]{Gueintegral}. 
Recall that a Verra fourfold $X$ is a smooth double cover
\begin{equation*}
\pi\colon X\longrightarrow\mathbb P^2\times\mathbb P^2
\end{equation*}
branched over a divisor of bidegree~$(2,2)$. It has two associated quadric surface fibrations which we denote $p_i=\operatorname{pr}_i\circ\pi:X\longrightarrow\mathbb P^2$ for $i=1,2$.

\begin{proposition}\label{prop:verra-lattice}
For a very general Verra fourfold $X$ one has
$$
T_X(1)\cong
U\oplus U(2)\oplus E_8(-1)^{\oplus2}\oplus\langle-2\rangle,
\quad 
D(T_X)\cong(\ZZ/2\ZZ)^3, \quad \ell(D(T_X))=3.
$$
$T_X$ is primitive in $\Lambda_X$ and has $\rk T_X = 21$. 
Furthermore $\End_{\mathrm{Hdg}}
\bigl(T_X\otimes_{\ZZ}\QQ\bigr)
= \QQ\cdot \id.$
\end{proposition}

\begin{proof}
 $T_X$ is primitive and nondegenerate by definition as $A_X\simeq U(2)\oplus\langle2\rangle$ is nondegenerate. The rank and discriminant group assertions follow from the stated lattice identification on the Hodge structure of $T_X$, which is provided in~\cite[Section~4]{Iliev_2017}. 
 
Let $S$ be the degree two K3 surface associated with either quadric surface fibration $p_i$. As rational Hodge structures we have 
$$
T_{X,\QQ}(1)\simeq T_{S,\QQ}\ceq T_S\otimes_{\ZZ}\QQ.
$$
As a consequence of~\cite[Proposition~2.10(ii)]{laterveer}, the associated degree two K3 surfaces vary dominantly in their moduli space. Hence $S$ is very general when $X$ is very general. For such an $S$, the Hodge group of $T(S)_{\QQ}$ is $\operatorname{SO}(T(S)_{\QQ})$; see \cite[Chapter~6, Corollary~4.7 and Theorem~4.9]{Huybrechts_2016}. Thus $\End_{\mathrm{Hdg}}(T(S)_{\QQ})=\QQ\cdot\operatorname{id}.$
\end{proof}

\begin{proof}[Proof of Corollary~\ref{cor:verra}]
By Proposition~\ref{prop:verra-lattice} a very general Verra fourfold $X$ satisfies \eqref{eq:conditions} with $\rk T_X=21$, $h^{3,1}(X)=1$ and $\ell(D(T_X))=3$. If $X$ were rational, Proposition~\ref{pro:quantum-integral}(2) would give $3\le\rho(S)\le10+12-21=1$.
\end{proof}

\begin{remark}Various works have compared K\"uchle fourfolds with Verra fourfolds. In particular, Benedetti--Fay--Gu\'er\'e--Manivel--Perrin proved that Hodge-general K\"uchle fourfolds of type \kt{c5} are not birational to any Verra fourfold~\cite[Corollary~5.3]{benedetti2026atomiccriterionirrationalityquantum} and the author proved the same for very general K\"uchle fourfolds of type \kt{c7}~\cite[Corollary~1.6]{fay2026hasseobstructionsrationalityspecial}. Benedetti--Gu\'er\'e--Manivel--Perrin studied the birational geometry of Verra fourfolds in~\cite{benedetti2026quantumcohomologybirationalgeometry}. For very general symmetric Verra fourfolds, equivariant irrationality with respect to the factor-swap involution was proved in~\cite[Theorem~6.4]{fay2026equivariantirrationalitygeneralsymmetric}. More recently, Jacovskis--Meinsma--Moschetti gave a construction of cubic fourfolds birational to Verra fourfolds in \cite{jacovskis2026verrainvolutionspecialcubic}. 
\end{remark}

\section{Code Availability}\label{sec:code}
The code and data associated with the Quantum Lefschetz calculations in Section~\ref{sec:comp} can be found at \url{https://github.com/aideenfay/kuchle-quantum-lefschetz}. Our SageMath implementation of Atiyah--Bott equivariant localization for the line counts in Table~\ref{tab:kuchle-line-counts} is available at \url{https://github.com/aideenfay/grassmannian-lines}. 
\section{Acknowledgments}
The author thanks Tom Coates for his assistance and many valuable discussions, and gratefully acknowledges EPSRC grant EP/Y028872/1. Parts of this work were assisted by GPT-5.5 and GPT-6 Astra in the form of finding gaps in proofs and checking the final paper for mistakes. 

\bibliographystyle{alpha}
\bibliography{refs}

@article{Inoue_2017,
   title={{I}-functions of {C}alabi–{Y}au {3}-folds in {G}rassmannians},
   volume={11},
   ISSN={1931-4531},
   url={http://dx.doi.org/10.4310/CNTP.2017.v11.n2.a2},
   DOI={10.4310/cntp.2017.v11.n2.a2},
   number={2},
   journal={Communications in Number Theory and Physics},
   publisher={International Press of Boston},
   author={Inoue, Daisuke and Ito, Atsushi and Miura, Makoto},
   year={2017},
   pages={273–309} }

@misc{golyshev2008spectrastrains,
      title={{S}pectra and strains}, 
      author={V. Golyshev},
      year={2008},
      eprint={0801.0432},
      archivePrefix={arXiv},
      primaryClass={hep-th},
      url={https://arxiv.org/abs/0801.0432}, 
}

@misc{Iritani2017QuantumLefschetz,
  author = {Iritani, Hiroshi},
  title  = {{Q}uantum ring and {Q}uantum {Lefschetz}},
  year   = {2017},
  url    = {https://www.math.kyoto-u.ac.jp/~iritani/quantum_ring_quantum_Lefschetz.pdf}
}

@article{Deligne1972,
author = {Deligne, Pierre},
journal = {Inventiones mathematicae},
pages = {206-226},
title = {La conjecture de {W}eil pour les surfaces {K3}.},
url = {http://eudml.org/doc/142129},
volume = {15},
year = {1971/72},
}

@misc{hu2025bigquantumcohomologyfano,
      title={Big quantum cohomology of {F}ano complete intersections}, 
      author={Xiaowen Hu},
      year={2025},
      eprint={1501.03683},
      archivePrefix={arXiv},
      primaryClass={math.AG},
      url={https://arxiv.org/abs/1501.03683}, 
}

@misc{fay2026equivariantirrationalitygeneralsymmetric,
      title={Equivariant irrationality of very general symmetric {V}erra fourfolds}, 
      author={Aideen Fay},
      year={2026},
      eprint={2605.30439},
      archivePrefix={arXiv},
      primaryClass={math.AG},
      url={https://arxiv.org/abs/2605.30439}, 
}

@article{laterveer,
  title={Algebraic cycles and {V}erra fourfolds},
  author={Robert Laterveer},
  journal={Tohoku Mathematical Journal, Second Series},
  volume={72},
  number={3},
  pages={451-485},
  year={2020},
  doi={10.2748/tmj/1601085625}
}

@article{Coates01062020,
	author = {Coates, Tom and Galkin, Sergey and Kasprzyk, Alexander and Strangeway, Andrew},
	doi = {10.1080/10586458.2018.1448018},
	eprint = {https://doi.org/10.1080/10586458.2018.1448018},
	journal = {Experimental Mathematics},
	number = {2},
	pages = {183--221},
	publisher = {Taylor \& Francis},
	title = {Quantum Periods for Certain Four-Dimensional {F}ano Manifolds},
	url = {https://doi.org/10.1080/10586458.2018.1448018},
	volume = {29},
	year = {2020}}

@article{deligne1980weil2,
     author = {Deligne, Pierre},
     title = {La conjecture de {Weil} : {II}},
     journal = {Publications Math\'ematiques de l'IH\'ES},
     pages = {137--252},
     year = {1980},
     publisher = {Institut des Hautes \'Etudes Scientifiques},
     volume = {52},
     doi = {10.1007/BF02684780},
     mrnumber = {601520},
     zbl = {0456.14014},
     language = {fr},
     url = {https://www.numdam.org/articles/10.1007/BF02684780/}
}

@incollection{katz2004larsen,
  author    = {Katz, Nicholas M.},
  title     = {Larsen's alternative, moments, and the monodromy of {L}efschetz pencils},
  booktitle = {Contributions to automorphic forms, geometry, and number theory},
  publisher = {Johns Hopkins University Press},
  address   = {Baltimore, MD},
  pages     = {521--560},
  year      = {2004}
}

@misc{jacovskis2026verrainvolutionspecialcubic,
      title={The {V}erra involution on special cubic fourfolds}, 
      author={Augustinas Jacovskis and Reinder Meinsma and Riccardo Moschetti},
      year={2026},
      eprint={2610.32952},
      archivePrefix={arXiv},
      primaryClass={math.AG},
      url={https://arxiv.org/abs/2610.32952}, 
}

@book{Huybrechts_2016, place={Cambridge}, series={Cambridge Studies in Advanced Mathematics}, title={Lectures on K3 Surfaces}, publisher={Cambridge University Press}, author={Huybrechts, Daniel}, year={2016}, collection={Cambridge Studies in Advanced Mathematics}}

@article{deligne,
	author = {Deligne, Pierre},
	date = {1974/12/01},
	doi = {10.1007/BF02684373},
	id = {Deligne1974},
	isbn = {1618-1913},
	journal = {Publications Math{\'e}matiques de l'Institut des Hautes {\'E}tudes Scientifiques},
	number = {1},
	pages = {273--307},
	title = {La conjecture de {W}eil. {I}},
	url = {https://doi.org/10.1007/BF02684373},
	volume = {43},
	year = {1974}}

@misc{AIF_0__0_0_A50_0,
     author = {Bernardara, Marcello and Fatighenti, Enrico and Manivel, Laurent and Tanturri, Fabio},
     title = {{F}ano fourfolds of {K3} type},
     journal = {Annales de l'Institut Fourier},
     year = {2026},
     publisher = {Association des Annales de l{\textquoteright}institut Fourier},
     doi = {10.5802/aif.3761},
     language = {en},
}

@misc{fay2026hasseobstructionsrationalityspecial,
      title={Hasse Obstructions to Rationality for Special Fourfolds}, 
      author={Aideen Fay},
      year={2026},
      eprint={2609.06759},
      archivePrefix={arXiv},
      primaryClass={math.AG},
      url={https://arxiv.org/abs/2609.06759}, 
}

@article{nikulin,
  title={Integral Symmetric Bilinear Forms and Some of Their Applications},
  author={Viacheslav V. Nikulin},
  journal={Mathematics of The Ussr-izvestiya},
  year={1980},
  volume={14},
  pages={103-167},
  url={https://api.semanticscholar.org/CorpusID:123044125}
}

@misc{meng2024quantumcohomologyfanoquiver,
      title={Quantum Cohomology of a {F}ano Quiver Moduli Space}, 
      author={Junyu Meng},
      year={2024},
      eprint={2412.15987},
      archivePrefix={arXiv},
      primaryClass={math.AG},
      url={https://arxiv.org/abs/2412.15987}, 
}

@misc{benedetti2026atomiccriterionirrationalityquantum,
      title={An atomic criterion for irrationality without quantum computations}, 
      author={Vladimiro Benedetti and Aideen Fay and Jérémy Guéré and Laurent Manivel and Nicolas Perrin},
      year={2026},
      eprint={2607.26718},
      archivePrefix={arXiv},
      primaryClass={math.AG},
      url={https://arxiv.org/abs/2607.26718}, 
}

@misc{alexeev2026atomicdecompositioncompleteintersection,
      title={On the atomic decomposition of complete intersection in flag varieties}, 
      author={Boris Alexeev and Leonardo F. Cavenaghi and Giovane Galindo and Bogdan Georgiev and Ludmil Katzarkov and Pedro Antonio Muniz Martins},
      year={2026},
      eprint={2609.21914},
      archivePrefix={arXiv},
      primaryClass={math.AG},
      url={https://arxiv.org/abs/2609.21914}, 
}

@article{ATIYAH19841,
	author = {M.F. Atiyah and R. Bott},
	doi = {https://doi.org/10.1016/0040-9383(84)90021-1},
	issn = {0040-9383},
	journal = {Topology},
	number = {1},
	pages = {1-28},
	title = {The moment map and equivariant cohomology},
	url = {https://www.sciencedirect.com/science/article/pii/0040938384900211},
	volume = {23},
	year = {1984}}

@article{Kuznetsov_2015,
   title={On {K}üchle varieties with {P}icard number greater than 1},
   volume={79},
   ISSN={1468-4810},
   url={http://dx.doi.org/10.1070/IM2015v079n04ABEH002758},
   DOI={10.1070/im2015v079n04abeh002758},
   number={4},
   journal={Izvestiya: Mathematics},
   publisher={Steklov Mathematical Institute},
   author={Kuznetsov, A G},
   year={2015},
   month=Aug, pages={698–709} }

@misc{benedetti2026quantumcohomologyirrationalitygushelmukai,
      title={Quantum cohomology and irrationality of {G}ushel--{M}ukai fourfolds}, 
      author={Vladimiro Benedetti and Laurent Manivel and Nicolas Perrin},
      year={2026},
      eprint={2603.17487},
      archivePrefix={arXiv},
      primaryClass={math.AG},
      url={https://arxiv.org/abs/2603.17487}, 
}

@misc{Gueirrational,
      title={On the irrationality of cubic fourfolds}, 
author={J{\'e}r{\'e}my Gu{\'e}r{\'e}},
year={2026},
      eprint={2603.04518},
      archivePrefix={arXiv},
      primaryClass={math.AG},
      url={https://arxiv.org/abs/2603.04518}, 
}

@article{Kuznetsov_2016,
   title={{K}\"uchle fivefolds of type c5},
   volume={284},
   ISSN={1432-1823},
   url={http://dx.doi.org/10.1007/s00209-016-1707-9},
   DOI={10.1007/s00209-016-1707-9},
   number={3-4},
   journal={Mathematische Zeitschrift},
   publisher={Springer Science and Business Media LLC},
   author={Kuznetsov, Alexander},
   year={2016},
 pages={1245–1278} }

@article{Lewis1985TheCH,
  title={The cylinder homomorphism associated to quintic fourfolds},
  author={James D. Lewis},
  journal={Compositio Mathematica},
  year={1985},
  volume={56},
  pages={315-329},
  url={https://api.semanticscholar.org/CorpusID:198461810}
}

@book{Fulton1984,
  author    = {Fulton, William},
  title     = {{I}ntersection {T}heory},
  series    = {Ergebnisse der Mathematik und ihrer Grenzgebiete. 3. Folge / A Series of Modern Surveys in Mathematics},
  volume    = {2},
  publisher = {Springer-Verlag},
  address   = {Berlin, Heidelberg},
  year      = {1984},
  doi       = {10.1007/978-3-662-02421-8},
  url       = {https://doi.org/10.1007/978-3-662-02421-8}
}

@article{Lew88,
 ISSN = {00029327, 10806377},
 URL = {http://www.jstor.org/stable/2374540},
 author = {James D. Lewis},
 journal = {American Journal of Mathematics},
 number = {1},
 pages = {77--114},
 publisher = {Johns Hopkins University Press},
 title = {The Cylinder Correspondence for Hypersurfaces of Degree $n$ in {$\PP^n$}},
 urldate = {2026-11-12},
 volume = {110},
 year = {1988}
}

@misc{li1998virtualmodulicyclesgromovwitten,
      title={Virtual moduli cycles and {G}romov--{W}itten invariants of algebraic varieties}, 
      author={Jun Li and Gang Tian},
      year={1998},
      eprint={alg-geom/9602007},
      archivePrefix={arXiv},
      primaryClass={alg-geom},
      url={https://arxiv.org/abs/alg-geom/9602007}, 
}

@article{Kontsevich_1994,
   title={{G}romov--{W}itten classes, quantum cohomology, and enumerative geometry},
   volume={164},
   ISSN={1432-0916},
   url={http://dx.doi.org/10.1007/BF02101490},
   DOI={10.1007/bf02101490},
   number={3},
   journal={Communications in {M}athematical {P}hysics},
   publisher={Springer Science and Business Media LLC},
   author={Kontsevich, M. and Manin, Yu.},
   year={1994},
   month=Aug, pages={525–562} }

@article{shimada,
	author = {Ichiro Shimada},
	doi = {10.2969/jmsj/04240719},
	journal = {Journal of the Mathematical Society of Japan},
	number = {4},
	pages = {719 -- 738},
	publisher = {Mathematical Society of Japan},
	title = {{On the cylinder homomorphisms of {F}ano complete intersections}},
	url = {https://doi.org/10.2969/jmsj/04240719},
	volume = {42},
	year = {1990}}

@article{Coates2001QuantumRL,
  title={Quantum {R}iemann–{R}och, {L}efschetz and {S}erre},
  author={Tom Coates and Alexander B. Givental},
  journal={Annals of Mathematics},
  year={2007},
  volume={165},
  pages={15-53},
  url={https://api.semanticscholar.org/CorpusID:12272063}
}

@article{Iliev_2017,
   title={Hyper-{K}\"ahler fourfolds and {K}ummer surfaces},
   volume={115},
   ISSN={0024-6115},
   url={http://dx.doi.org/10.1112/plms.12063},
   DOI={10.1112/plms.12063},
   number={6},
   journal={Proceedings of the London Mathematical Society},
   publisher={Wiley},
   author={Iliev, Atanas and Kapustka, Grzegorz and Kapustka, Michał and Ranestad, Kristian},
   year={2017},
   month=sep, pages={1276–1316} }

@misc{benedetti2026quantumcohomologybirationalgeometry,
      title={Quantum cohomology and birational geometry of {V}erra fourfolds}, 
      author={Vladimiro Benedetti and Jérémy Guéré and Laurent Manivel and Nicolas Perrin},
      year={2026},
      eprint={2605.30450},
      archivePrefix={arXiv},
      primaryClass={math.AG},
      url={https://arxiv.org/abs/2605.30450}, 
}

@article{Collino_1999,
   title={On the Structure of the Small Quantum Cohomology Rings of Projective Hypersurfaces},
   volume={206},
   ISSN={1432-0916},
   url={http://dx.doi.org/10.1007/s002200050701},
   DOI={10.1007/s002200050701},
   number={1},
   journal={Communications in Mathematical Physics},
   publisher={Springer Science and Business Media LLC},
   author={Collino, Alberto and Jinzenji, Masao},
   year={1999},
   month=sep, pages={157–183} }

@misc{katzarkov2026birationalinvariantshodgestructures,
      title={Birational Invariants from {H}odge Structures and Quantum Multiplication}, 
      author={Ludmil Katzarkov and Maxim Kontsevich and Tony Pantev and Tony Yue YU},
      year={2026},
      eprint={2508.05105},
      archivePrefix={arXiv},
      primaryClass={math.AG},
      url={https://arxiv.org/abs/2508.05105}, 
}

@misc{Gueintegral,
      title={Integral {H}odge structures and atoms of fourfolds}, 
author={J{\'e}r{\'e}my Gu{\'e}r{\'e}},
year={2026},
      eprint={2610.40811},
      archivePrefix={arXiv},
      primaryClass={math.AG},
      url={https://arxiv.org/abs/2610.40811}, 
}

@article{kuchle,
	author = {K{\"u}chle, Oliver},
	date = {1995/01/01},
	doi = {10.1007/BF02571923},
	id = {K{\"u}chle1995},
	isbn = {1432-1823},
	journal = {Mathematische Zeitschrift},
	number = {1},
	pages = {563--575},
	title = {On {F}ano 4-folds of index 1 and homogeneous vector bundles over {G}rassmannians},
	url = {https://doi.org/10.1007/BF02571923},
	volume = {218},
	year = {1995}}

@misc{coates2014quantumlefschetzprinciplevector,
      title={The {Q}uantum {L}efschetz Principle for Vector Bundles as a Map Between {G}ivental Cones}, 
      author={Tom Coates},
      year={2014},
      eprint={1405.2893},
      archivePrefix={arXiv},
      primaryClass={math.AG},
      url={https://arxiv.org/abs/1405.2893}, 
}
\clearpage
\appendix
\section{Line counts by Atiyah--Bott Equivariant Localization}
\label{app:lines}
Let $X\subset G=\Gr(k,n)$ be a fourfold from K\"uchle's list~\cite{kuchle} defined by the bundle $\mathcal{F}\oplus\cO(1)$, and let $Y'$ be the fivefold cut out by $\mathcal{F}$. Thus $X$ is a hyperplane
section of $Y'$ in the Pl\"ucker embedding $G\hookrightarrow\PP(\bigwedge^k\mathbb C^n)$. For a general point $p\in Y'$, every line $\ell \in G$ through $p$ is parameterized by $\PP^{k-1}\times\PP^{n-k-1}$, and those contained in $Y'$ are the zeros of a section of the bundle $\mathcal E$ with fibres $$\mathcal E_\ell = H^0\bigl(\ell,\mathcal{F}|_\ell(-p)\bigr).$$
Since $\mathcal E$ has rank $c_1(F)\cdot\ell=n-2$, the number of lines through $p$ is
\[
  a_X =
  \int_{\PP^{k-1}\times\PP^{n-k-1}} c_{n-2}(\mathcal E).
\]
We compute this integral using Atiyah--Bott equivariant localization~\cite[Equation~3.8]{ATIYAH19841} and record the results in Table~\ref{tab:kuchle-line-counts}. See Section~\ref{sec:code} for information on code availability. 
\begin{table}[htb!]
  \centering
  \renewcommand{\arraystretch}{1.15}
  \setlength{\tabcolsep}{6pt}
  \begin{tabular}{@{}lclc@{}}
    \toprule
    K\"uchle type & G & $\mathcal F$ & $a_X$ \\
    \midrule
    \kt{b1}  & $\Gr(2,5)$
      & $\mathcal O(3)$ & $18$ \\
    \kt{b5}  & $\Gr(2,6)$
      & $\mathcal S^\vee(1)\oplus\mathcal O(1)$ & $7$ \\
    \kt{b6}  & $\Gr(2,6)$
      & $\mathcal O(1)^{\oplus2}\oplus\mathcal O(2)$ & $8$ \\
    \kt{b7}  & $\Gr(2,7)$
      & $\mathcal O(1)^{\oplus5}$ & $5$ \\
    \kt{b8}  & $\Gr(2,7)$
      & $\operatorname{Sym}^2\mathcal S^\vee
         \oplus\mathcal O(1)^{\oplus2}$ & $4$ \\
    \kt{b11} & $\Gr(2,8)$
      & $\mathcal Q^\vee(1)\oplus\mathcal O(1)$ & $3$ \\
    \midrule
    \kt{c1}  & $\Gr(3,6)$
      & $\mathcal O(1)^{\oplus4}$ & $6$ \\
    \kt{c2}  & $\Gr(3,6)$
      & $\bigwedge^2\mathcal S^\vee\oplus\mathcal O(2)$ & $8$ \\
    \kt{c4}  & $\Gr(3,7)$
      & $\operatorname{Sym}^2\mathcal S^\vee
         \oplus\mathcal O(1)$ & $0$ \\
    \kt{c5}  & $\Gr(3,7)$
      & $\bigwedge^2\mathcal S^\vee
         \oplus\mathcal Q^\vee(1)$ & $4$ \\
    \kt{c6}  & $\Gr(3,7)$
      & $\bigl(\bigwedge^2\mathcal S^\vee\bigr)^{\oplus2}
         \oplus\mathcal O(1)$ & $5$ \\
    \kt{c7}  & $\Gr(3,8)$
      & $\bigwedge^3\mathcal Q$ & $0$ \\
    \midrule
    \kt{d1}  & $\Gr(4,8)$
      & $\operatorname{Sym}^2\mathcal S^\vee
         \oplus\mathcal O(1)$ & $0$ \\
    \kt{d3}  & $\Gr(5,10)$
      & $\bigl(\bigwedge^2\mathcal S^\vee\bigr)^{\oplus2}$ & $5$ \\
    \bottomrule
  \end{tabular}
  \caption{The numbers $a_X$ of lines through a general point of the
  fivefold $Y'$ defined by $\mathcal{F}$, containing the fourfold $X$ as a hyperplane section. }
  \label{tab:kuchle-line-counts}
\end{table}

\end{document}